\documentclass[letter,11pt]{amsart}
\usepackage{geometry,wasysym}
\usepackage[utf8]{inputenc}
\usepackage{mathpazo}  
\usepackage[english]{babel}
\usepackage[usenames,dvipsnames]{color}
\usepackage{url}
\usepackage{pdflscape}
\usepackage{amsmath,amssymb,latexsym,mathrsfs,amssymb,amsthm}
\usepackage{enumerate}
\usepackage{mathtools} 
\usepackage{epstopdf}
\usepackage{graphicx}
\usepackage{cite}
\usepackage[all]{xy}
\usepackage{setspace}
\usepackage[normalem]{ulem}
\usepackage[norelsize]{algorithm2e}
\usepackage{newfloat}
\DeclareFloatingEnvironment[
    fileext=los,
    listname=List of Schemes,
    name=Notes on Algorithm,
    placement=tbhp,
]{scheme}

\newcommand\cmd[1]{\setbox0=\hbox{x}\indent\texttt{ >> #1}}

\usepackage[colorlinks=true,linkcolor=PineGreen, citecolor=cyan, urlcolor=PineGreen, pagebackref]{hyperref}

\usepackage{color,graphicx,tikz, tkz-euclide, float, pgfplots, overpic}
\usepackage{subfigure}
\usetikzlibrary{arrows,calc, decorations.markings, positioning, fadings}
\usetkzobj{all}

\newtheorem{theorem}{Theorem}[section]

\theoremstyle{definition}

\theoremstyle{remark}
\newtheorem{remark}[theorem]{Remark}

\DeclareMathOperator{\sech}{sech}

\newcommand{\D}{\ensuremath{\,\mathrm{d}}}			
\newcommand{\defeq}{\vcentcolon=}

\makeatletter
\renewcommand*\env@matrix[1][\arraystretch]{%
  \edef\arraystretch{#1}%
  \hskip -\arraycolsep
  \let\@ifnextchar\new@ifnextchar
  \array{*\c@MaxMatrixCols c}}
\makeatother

\let\originalleft\left
\let\originalright\right
\renewcommand{\left}{\mathopen{}\mathclose\bgroup\originalleft}
\renewcommand{\right}{\aftergroup\egroup\originalright}

\title{Benchmarking Numerical Methods for Lattice Equations With the Toda Lattice}

\author[D.~Bilman]{Deniz Bilman}
\address{Deniz Bilman\\
     Department of Mathematics\\
     University of Michigan\\
     530 Church Street\\
     Ann Arbor, MI 48109
     }
\email{bilman@umich.edu}
\author[T.~Trogdon]{Thomas Trogdon}
\address{Thomas Trogdon\\
     University of California, Irvine\\
     Rowland Hall\\   
     Irvine, CA 92697
   }
\email{ttrogdon@uci.edu}
\date{\today}
\begin{document}

\begin{abstract}
We compare performances of well-known numerical time-stepping methods that are widely used to compute solutions of the doubly-infinite Fermi-Pasta-Ulam-Tsingou (FPUT) lattice equations. The methods are benchmarked according to (1) their accuracy in capturing the soliton peaks and (2) in capturing highly-oscillatory parts of the solutions of the Toda lattice resulting from a variety of initial data. The numerical inverse scattering transform method is used to compute a reference solution with high accuracy.  We find that benchmarking a numerical method on pure-soliton initial data can lead one to overestimate the accuracy of the method.
\end{abstract}

\maketitle

\section{Introduction}
\begin{figure}[htp]
\begin{tikzpicture}[scale=1.2]
\filldraw[CadetBlue] (-4,0) circle (2.2pt);
\filldraw[CadetBlue] (-3,0) circle (2.2pt);
\filldraw[CadetBlue] (-2,0) circle (2.2pt);
\filldraw[CadetBlue] (-1,0) circle (2.2pt);
\filldraw[CadetBlue] (0,0) circle (2.2pt);
\filldraw[CadetBlue] (1,0) circle (2.2pt);
\filldraw[CadetBlue] (2,0) circle (2.2pt);
\filldraw[CadetBlue] (3,0) circle (2.2pt);
\filldraw[CadetBlue] (4,0) circle (2.2pt);
\draw[thick,PineGreen,decoration={aspect=0.4, segment length=1mm, amplitude=1mm,coil},decorate,opacity=1,path fading=west] (-4.8,0)--(-4.2,0);
\draw[thick,PineGreen,opacity=1](-4.2,0)--(-4.08,0);

\draw[thick,PineGreen](-3.92,0)--(-3.8,0);
\draw[thick,PineGreen,decoration={aspect=0.4, segment length=1mm, amplitude=1mm,coil},decorate] (-3.8,0)--(-3.2,0);
\draw[thick,PineGreen](-3.2,0)--(-3.08,0);

\draw[thick,PineGreen](-2.92,0)--(-2.8,0);
\draw[thick,PineGreen,decoration={aspect=0.4, segment length=1mm, amplitude=1mm,coil},decorate] (-2.8,0)--(-2.2,0);
\draw[thick,PineGreen](-2.2,0)--(-2.08,0);

\draw[thick,PineGreen](-1.92,0)--(-1.8,0);
\draw[thick,PineGreen,decoration={aspect=0.4, segment length=1mm, amplitude=1mm,coil},decorate] (-1.8,0)--(-1.2,0);
\draw[thick,PineGreen](-1.2,0)--(-1.08,0);

\draw[thick,PineGreen](-0.92,0)--(-0.8,0);
\draw[thick,PineGreen,decoration={aspect=0.4, segment length=1mm, amplitude=1mm,coil},decorate] (-0.8,0)--(-0.2,0);
\draw[thick,PineGreen](-0.2,0)--(-0.08,0);

\draw[thick,PineGreen](0.08,0)--(0.2,0);
\draw[thick,PineGreen,decoration={aspect=0.4, segment length=1mm, amplitude=1mm,coil},decorate] (0.2,0)--(0.8,0);
\draw[thick,PineGreen](0.8,0)--(0.92,0);

\draw[thick,PineGreen](1.08,0)--(1.2,0);
\draw[thick,PineGreen,decoration={aspect=0.4, segment length=1mm, amplitude=1mm,coil},decorate] (1.2,0)--(1.8,0);
\draw[thick,PineGreen](1.8,0)--(1.92,0);

\draw[thick,PineGreen](2.08,0)--(2.2,0);
\draw[thick,PineGreen,decoration={aspect=0.4, segment length=1mm, amplitude=1mm,coil},decorate] (2.2,0)--(2.8,0);
\draw[thick,PineGreen](2.8,0)--(2.92,0);

\draw[thick,PineGreen](3.08,0)--(3.2,0);
\draw[thick,PineGreen,decoration={aspect=0.4, segment length=1mm, amplitude=1mm,coil},decorate,] (3.2,0)--(3.8,0);
\draw[thick,PineGreen](3.8,0)--(3.92,0);

\draw[thick,PineGreen,opacity=1](4.08,0)--(4.2,0);
\draw[thick,PineGreen,decoration={aspect=0.4, segment length=1mm, amplitude=1mm,coil},decorate,opacity=1,path fading=east] (4.2,0)--(4.8,0);

\node (x) at (0,-0.8){\small \hspace{2cm}$q_{n}(t)$: displacement of the $n^{\text{th}}$ particle};
\node (x2) at (0, -1.1){\small \hspace{2.6cm}from its equilibrium position};
\node (n) at (0,0.48){\small$n$};
\node (n1) at (1,0.5){\small $n+1$};
\node (n2) at (2,0.5){$\cdots$};
\node (n3) at (-1,0.5){\small$n-1$};
\node (n4) at (-2, 0.5){$\cdots$};
\draw (0,-0.2)--(0,-0.6);
\draw[->](0.1,-0.6)--(0.5,-0.6);
\draw[->](-0.1,-0.6)--(-0.5,-0.6);
\draw[black,thin,<->] plot[smooth] coordinates {
    (-0.6,-0.2)
     (-1,-0.45)
     (-1.67,-0.5)
   };
\node[black] (p) at (-3,-0.5){\small nonlinear springs};
\end{tikzpicture}
\caption{One-dimensional chain of particles with nearest neighbor interactions.}
\label{F:lattice}
\end{figure}
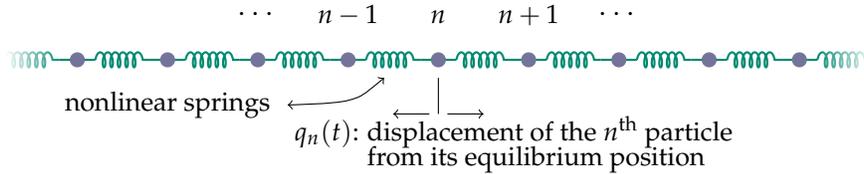
Consider the classical problem of one-dimensional infinite chain of particles on a line with nearest-neighbor interactions as depicted in Figure~\ref{F:lattice}. Assume that each particle has unit mass, and that there are no impurities, \emph{i.e}.\ the potential energies of the springs between the particles are identical. We let $V\colon \mathbb{R}\to \mathbb{R}$ denote the interaction potential between the neighboring particles. With these assumptions, the equations of motion that govern this particle system are given by Newton's Second Law of Motion:
\begin{equation}
\frac{\D^2}{\D t^2} q_n=V'(q_{n+1}-q_n) - V'(q_{n}-q_{n-1}),\quad n\in\mathbb{Z},
\label{eq:eom}
\end{equation}
where $q_n$ stands for the displacement of the $n^\mathrm{th}$ particle from its equilibrium position. Denoting by $p_n$ the momentum of the $n^{\text{th}}$ particle, \eqref{eq:eom} is equivalent to the system of first order differential equations:
\begin{equation}
\frac{\D p_n}{d t} =V{\left(q_{n+1} - q_{n}\right)})-V{\left(q_n - q_{n-1}\right)},\quad \frac{\D q_n}{d t} = p_n,\quad n\in\mathbb{Z}\,.
\label{eq:eom-pq}
\end{equation}
With the assumptions that $q_{n+1}-q_n \to 0$ and $p_n \to 0$ sufficiently fast as $|n|\to \infty$ (i.e.\ no motion at infinity), \eqref{eq:eom-pq} is a Hamiltonian system of equations
\begin{equation}
\label{eq:Hamiltonian}
\frac{d p_n}{d t} = - \frac{\partial \mathcal{H}(p,q)}{\partial q_n},\quad
\frac{d q_n}{d t} = \frac{\partial \mathcal{H}(p,q)}{\partial p_n},\quad n\in\mathbb{Z}\,,
\end{equation}
with the Hamiltonian functional $\mathcal{H}(p,q)$:
\begin{equation}
\mathcal{H}(p,q)=\sum_{n\in\mathbb{Z}}\tfrac{1}{2}p_n^2 + V\left(q_{n+1}-q_n\right).
\label{eq:Toda-ham}
\end{equation}
Such nearest-neighbor interacting particle systems with nonlinear interaction forces (anharmonic potentials) include the systems studied in the famous experiment at Los Alamos by Fermi, Pasta, Ulam, and Tsingou \cite{Fermi1955a} in 1953, which a decade later led to discovery of solitons by Zabusky and Kruskal \cite{Zabusky1965a}. Such lattice equations model various physical phenomena with a multitude of applications \cite{Sato2003}.  From a purely mathematical perspective these lattices are used to investigate Poincar\'e recurrence, chaos, and nonlinear wave phenomena (interaction of solitary waves, solitary wave resolution, see \cite{Bilman2016} and the references therein).  The FPUT-type systems \eqref{eq:eom-pq} are still an active area of mathematical research \cite{Vainchtein2016,Hoffman2017,Gaison2014}.

The particular choice $V(r)=V_{\text{Toda}}(r)\defeq e^{-r} + r - 1$ in \eqref{eq:eom-pq} results in an infinite dimensional, continuous time - discrete space completely integrable system: the celebrated Toda lattice \cite{Toda1967a}.  The main purpose of this work is to compare the performances of various well known numerical time-stepping methods that are widely used to compute the solution of the Cauchy initial value problem for \eqref{eq:eom}. The accuracy of a time-stepping method is most easily inspected when there are exact solutions available at hand, which is of course, rarely the case. Often, numerical analysts make use of nonlinear integrable wave models that possess classes of explicit solutions (\emph{e.g.} the Toda lattice, the Korteweg-deVries (KdV) equation or the nonlinear Schr\"odinger (NLS) equation).  Then one has the luxury of being able to test their numerical scheme against exact formulae. Many of these infinite-dimensional integrable systems feature dispersive radiation and/or oscillatory tails.  This is a critically important\footnote{Generic initial data for the Toda lattice gives rise to dispersive radiation.} oscillatory component of the solution that decays slowly to the background as $t \to \infty$. It is often the case that exact solutions formulae are for solitary waves such as breathers or solitons --- coherent structures that are localized in space without oscillatory tails or dispersive radiation. In this paper, we set out to investigate questions such as: Which time-stepping methods capture the solitons with more accuracy? How do they perform when computing highly-oscillatory solutions? How does their performance depend on the solution itself? Our strategy is as follows. Using a numerical inverse scattering transform (IST)  method \cite{Bilman2017}, we can accurately construct solutions of the Toda lattice for each $(n,t)$ without any time stepping, for arbitrarily large values of $t$, with high accuracy.  Indeed, one can expect to maintain relative accuracy for large $t$ \cite{TrogdonSONNSD}. This opens the door to benchmarking time-stepping methods on solutions of different characters (oscillatory, has solitons, no solitons, etc.), without restricting oneself to those with exact formulae.  Furthermore, it is reasonable to expect behavior of a numerical method for the Toda lattice to extend to other FPUT-type systems.

{From a convergence (as time-step size tends to $0$) and computational complexity point of views, comparison of numerical time-stepping methods is a well-trodden path, see for example, ~\cite{Hull1972} for a detailed study in this direction, or the more recent survey article \cite{Butcher2000} by Butcher and the references therein.}  Using the numerical IST method, in this work we aim to add a new dimension to such studies.

The paper is organized as follows. In the remainder of this section we discuss specifics of the Toda lattice.  In Section~\ref{s:num-data} we summarize the numerical IST procedure, give an overview of the numerical time-stepping methods used in this work, and present the explicit formulae for the initial data used in the numerical experiments. In Section~\ref{s:2nd-order}, we compare performances of second order methods. In Section~\ref{s:higher-order}, we compare performances of higher order methods.  We then draw some conclusions.

\subsection{Properties of the Toda lattice}

Complete integrability of the Toda lattice was proven by H.~Flaschka \cite{Flaschka1974} and S.~V.~Manakov \cite{Manakov1975} in 1974, independently and simultaneously, by realizing that the system possesses a Lax pair and \eqref{eq:eom-pq} with $V(r) =V_\mathrm{T}(r)$ is in one-to-one correspondence with isospectral deformations of Jacobi (symmetric, tridiagonal with positive off-diagonal elements) matrices. Indeed, through the bijection
\begin{equation}
a_n \defeq\frac{1}{2}e^{-(q_{n+1} - q_{n})/2},\quad b_n \defeq -\frac{1}{2}p_n\,,
\label{eq:ab}
\end{equation}
the equations of motion \eqref{eq:eom} Toda lattice take the form:
\begin{equation}
\frac{\D }{\D t} a_n = a_n \left(b_{n+1} - b_n \right),\quad
\frac{\D }{\D t} b_n  = 2 \left(a_n^2 - a_{n-1}^2 \right),\quad n\in\mathbb{Z}.
\label{eq:eom-ab}
\end{equation}
Defining the following second order linear difference operators $\mathbf{L}$ and $\mathbf{P}$
on the Hilbert space $\ell^{2}(\mathbb{Z})$ of square-summable sequences:
\begin{equation}
\begin{aligned}
(\mathbf{L}\phi)_n &\defeq a_{n-1}\phi_{n-1} + b_n \phi_n + a_n \phi_{n+1}\\
(\mathbf{P}\phi)_n &\defeq -a_{n-1}\phi_{n-1} + a_n \phi_{n+1}\,, 
\end{aligned}
\label{eq:PL}
\end{equation}
it can be verified that \eqref{eq:eom} is equivalent to the Lax equation \cite{Lax1968a}
\begin{equation}
\frac{\D }{\D t} \mathbf{L} = [\mathbf{P},\mathbf{L}]\defeq \mathbf{P}\mathbf{L}-\mathbf{L}\mathbf{P}.
\label{eq:lax}
\end{equation}
The operators $(\mathbf{P},\mathbf{L})$ are called a Lax pair and in the standard basis, $\mathbf{L}$ is a doubly-infinite Jacobi matrix. The Lax pair \eqref{eq:lax} constitutes the basis of to the IST method to solve the Cauchy initial value problem for \eqref{eq:eom-ab} for sufficiently decaying initial data (see, for example, \cite{Bilman2017a} for a recent survey of the IST for the Toda lattice).    

A numerical IST method was recently developed by the authors \cite{Bilman2017} for the doubly-infinite Toda lattice.  Implementations  for other integrable systems can be found  in \cite{TrogdonSONNSD,TrogdonSONLS,TrogdonSOKdV,TrogdonFiniteGenus} and these are summarized in \cite{TrogdonSOBook}.  An implementation for the numerical IST method can be found at \cite{ISTPackage}.  The method works, loosely speaking, by performing the following steps:
\begin{enumerate}
\item Compute the spectral data:  This involves solving the eigenvalue problem $\mathbf L \varphi = \lambda \varphi$ for for bounded eigenfunctions $\varphi$.  Assuming $(a_n,b_n) \to (1/2,0)$ at an appropriately rapid rate as $|n|\to\infty$, the spectrum consists of the interval $[-1,1]$ and a finite number of simple eigenvalues in $\mathbb R \setminus [-1,1]$.  One also computes, a function $R(z)$, defined on the spectrum, which is directly related to the spectral measure for $\mathbf L$.
\item Solve the inverse problem: Once $R(z)$ is known, for each $(n,t)$ there exists a contour $\Gamma = \Gamma(n,t) \subset \mathbb C$, a function $G(z;n,t)$, $G: \Gamma \to \mathbb C^{2\times 2}$ and an integral equation
  \begin{align*}
    \left(I - \frac{G}{2}\right) U + \mathcal H_{\Gamma}(GU) = G - I, \quad U: \Gamma \to \mathbb C^{2\times 2}.
  \end{align*}
  Here $\mathcal H_{\Gamma}$ is the Hilbert transform over $\Gamma$ and $I$ is the identity matrix.  This integral equation is solved for $U: \Gamma \to \mathbb C^{2\times 2}$ using the framework of S. Olver \cite{SORHFramework} (see \cite{RHPackage} for an implementation)   and the solution of the Toda lattice can be obtained in terms of integrals of $U$.  So, one can compute the map $(n,t) \mapsto U$ and hence $(n,t) \to (a_n(t),b_n(t))$ without time stepping.
\end{enumerate}

\section{Time-Stepping Methods}\label{s:num-data}
In this section we describe the time-stepping methods used to produce the results shown in this paper. The first three methods described below are classical second-order methods and the following methods are higher-order methods. All of the methods are used to numerically compute the solution $y(T)\in\mathbb{R}^N$ of the Cauchy initial value problem for an autonomous differential equation
\begin{equation}
\frac{dy}{dt} = f(y(t))\,,
\label{eq:model-ode}
\end{equation}
with initial condition $y(0)=y^0\in\mathbb{R}^N$ at a final time $t=T$.  We focus on classical methods to illustrate how one can use the numerical IST method to benchmark them.

Depending the nature of the scheme and the variables used, $f$ is an explicit, known function from $\mathbb{R}^N$ to $\mathbb{R}^N$ whose definition may differ among different methods (for example, in St\"ormer-Verlet method). We truncate the spatial domain (the lattice $\mathbb{Z}$) to $\mathbb{S}=\{-K,\dots,0,\dots,K\}$, with the appropriate boundary conditions for the variables $(a,b)$ in \eqref{eq:eom-ab} or $(p,q)$ in \eqref{eq:eom-pq} and $K>0$ chosen large so that the $y(T)$ does not feel the effect of the boundary. Thus, $N=2K+1$. In the descriptions that follow, $h>0$ denotes the time-step size. We use superscripts to denote the numerical iterates to avoid confusion with the indices of sequences: Given a computed solution $y^k$ at a time $t$, we denote the computed solution at time $t+h$ by $y^{k+1}$. $y^k_n$ is a scalar, the $n$-th element of the sequence $y^k \in \mathbb{R}^N$.

\subsection{Second-order Methods}
\subsubsection{Midpoint} We denote this second-order explicit method by \texttt{midpoint}. Given $y^k$, the algorithm to compute $y^{k+1}$ is as follows:\newline
\cmd{$y^{k+1} = y^{k} + h f\left(y^{k} + \tfrac{1}{2}h f\left(y^{k}\right) \right )$}\newline
Unless otherwise stated, we integrate \eqref{eq:eom-ab}. When we integrate \eqref{eq:eom-pq} instead we use the label\\ \texttt{midpointqp}.
\subsubsection{Second-order St\"ormer-Verlet} We denote this method by \texttt{sv2symp}. The St\"ormer-Verlet method is symplectic: It preserves the Hamiltonian \eqref{eq:Toda-ham} under exact arithmetic. We use the equations of motions \eqref{eq:eom} and set
\begin{equation*}
f_p(q)\defeq \left( e^{-\left(q_n - q_{n-1}\right)}-e^{-\left(q_{n+1} - q_{n}\right)} \right)_{n\in\mathbb{Z}} \,,\quad f_q(p)\defeq ( p_n)_{n\in\mathbb{Z}}
\end{equation*}
to denote the right hand sides. Since the Hamiltonian for the Toda lattice is separable, the method becomes explicit. Given $(p^k,q^k)$, the algorithm to compute $(p^{k+1},q^{k+1})$ is as follows:\newline
\cmd{$p^{k+1/2} = p^k + \tfrac{1}{2} h f_p(q^k)$}\\
\cmd{$q^{k+1} = q^k + h f_q\left(p^{k+1}\right)$}\\
\cmd{$p^{k+1} = p^{k+1} + \tfrac{1}{2} h f_p\left(q^{k+1}\right)$}

\subsection{Fourth-order Methods}

\subsubsection{Fourth-order Adams-Bashforth} We denote this method by \texttt{ab4}. The Adams-Bashforth method is a linear explicit multi-step method. To compute the first three iterates we perform lower-step Adams-Bashforth methods successively:\newline
\cmd{$y^0 = y^0$}\\
\cmd{$y^1 = y^0 + h f(y^0)$}\\
\cmd{$y^2 = y^1 + \tfrac{1}{2}h\left(- f(y^0) + 3 f(y^1)\right)$}\\
\cmd{$y^3 = y^2 + \tfrac{1}{12}h\left(5 f(y^0) - 16 f(y^1) + 23 f(y^2)\right)$}\newline
Then, given computed solutions $y^k, y^{k-1},y^{k-2},y^{k-3}$, $k\geq 3$ the algorithm to compute $y^{k+1}$ is as follows:\newline
\cmd{$y^{k+1} = y^{k} + \tfrac{1}{24}h\left(-9 f(y^{k-3}) + 37 f(y^{k-2}) - 59 f(y^{k-1}) + 55 f(y^{k})\right)$}

\subsubsection{Fourth-order Runge-Kutta} We denote this method by \texttt{rk4}. The Runge-Kutta method is a an explicit single-step method. Given the computed solution $y^k$, the algorithm \cite{Kutta1901,Runge1895} to compute $y^{k+1}$ is as follows:\newline
\cmd{$s_1 = h f(y^k)$}\\
\cmd{$s_2 = h f\left( y^0 + \tfrac{1}{2}s_1 \right)$}\\
\cmd{$s_3 = h f\left( y^0 + \tfrac{1}{2}s_2 \right)$}\\
\cmd{$s_4 = h f\left( y^0 + s_3 \right)$}\\
\cmd{$y^{k+1} = y^{k} + \tfrac{1}{6}h \left(s_1 + 2 s_2+ 2 s_3 + s_4\right)$}\\
Unless otherwise stated, we integrate \eqref{eq:eom-ab}. When we integrate \eqref{eq:eom-pq} we use the label \texttt{rk4qp}.
\subsubsection{Four-Five-order Runge-Kutta-Fehlberg} We denote this method by \texttt{rkf45}. This method is fourth order method but with only a single extra computation at each step the local error can be controlled by a fifth order method. This feature \cite{Fehlberg1970} is very practical for implementing adaptive-step-size methods. We do not implement any adaptive methods. Given the computed solution $y^k$, the algorithm \cite{Fehlberg1970, HairerBook} to compute $y^{k+1}$ is as follows:\newline
\cmd{$s_1= h f(y)$}\\
\cmd{$s_2 = f(y + h b_{21} s_1)$}\\
\cmd{$s_3 = f\left(y + h \left(b_{31} s_1 + b_{32} s_2 \right)\right)$}\\
\cmd{$s_4 = f\left(y + h \left(b_{41} s_1 + b_{42} s_2 +b_{43} s_3 \right)\right)$}\\
\cmd{$s_5 = f\left(y + h \left(b_{51} s_1 + b_{52} s_2 +b_{53} s_3 +b_{54} s_4 \right)\right)$}\\
\cmd{$s_6 = f\left(y + h \left(b_{61} s_1 + b_{62} s_2 +b_{63} s_3 +b_{64} s_4 + b_{65} s_5 \right)\right)$}\\
\cmd{$y^{k+1}= y^k + h\left(c_1 s_1 +  c_2 s_2 +c_3 s_3 + c_4 s_4 + c_5 s_5 + c_6 s_6\right)$}\newline
For the constants $b_{ij}$ and $c_j$, see the Butcher tableau provided in the Appendix~\ref{s:Butcher}. See \cite{Butcher1996} for a survey article by Butcher on Runge-Kutta methods.

\begin{remark}
For a detailed study of these methods and more, the reader may consult to textbooks \cite{HairerBook} or \cite{leveque}, for example.
\end{remark}

\section{Initial data considered in numerical experiments}
We now present the types of initial data (ID) used in the numerical experiments underlying the results presented in this paper.
\subsubsection{Purely dispersive ({\tt NoS})} 
\begin{figure}[h]
  \centering
  \includegraphics[width=\linewidth]{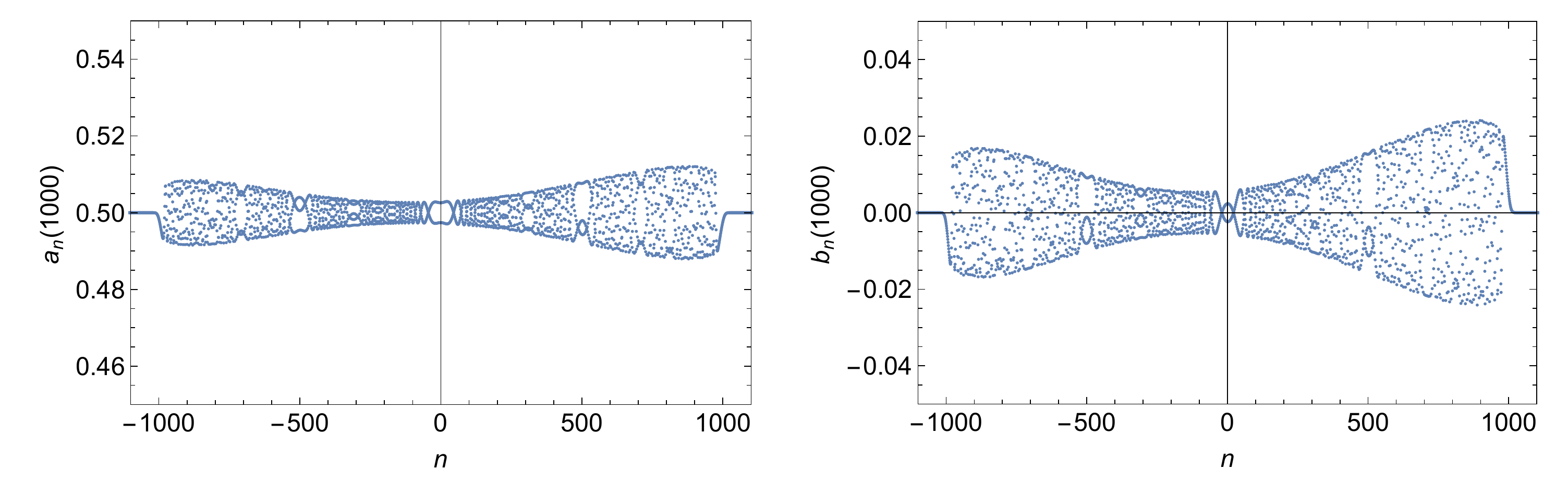}
  \caption{ The solution of the Toda lattice with \texttt{NoS} ID at $t = 1000$.}\label{f:NoS}
\end{figure}

We chose ID that gives rise to pure radiation: There are no solitons since the discrete spectrum of $\mathbf{L}$ is empty. The ID is given by
\begin{equation}
a_n = \frac{1}{2} - \frac{1}{4} e^{-n^2},\quad
b_n = \frac{1}{10}\sech(n),\quad n\in\mathbb{Z}.
\label{eq:datadisp}
\end{equation}
The solution at $t = 1000$ is shown in Figure~\ref{f:NoS}.

\subsubsection{Pure 1-soliton ({\tt PureS})} 
\begin{figure}[h]
  \centering
  \includegraphics[width=\linewidth]{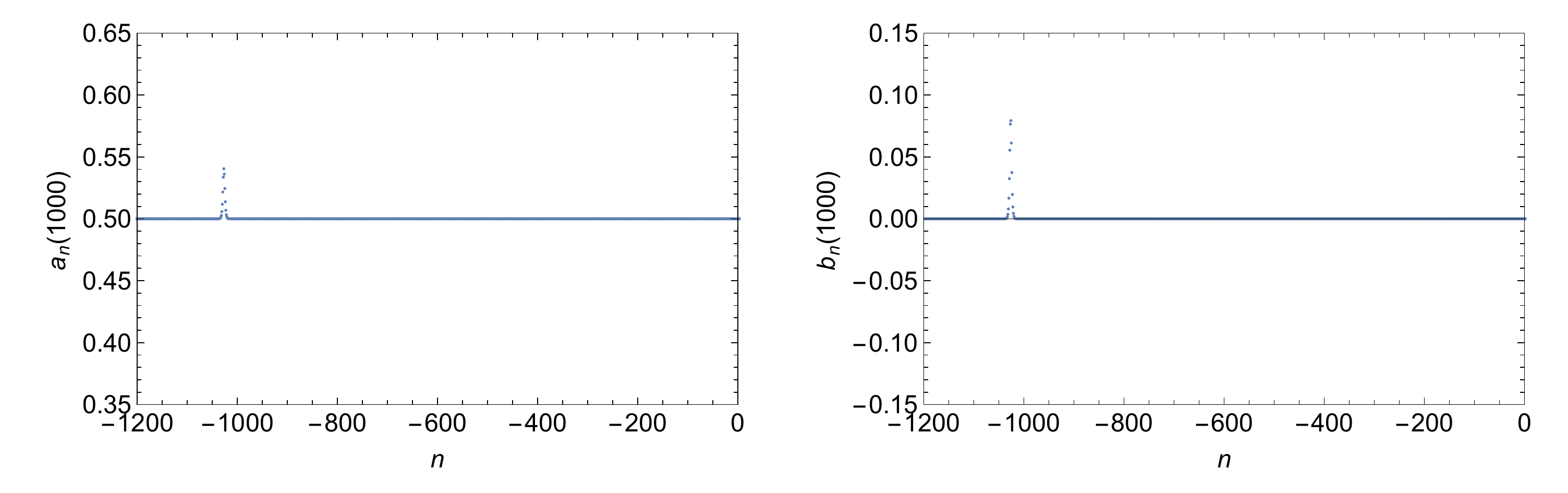}
  \caption{ The solution of the Toda lattice with {\tt PureS} ID at $t = 1000$.}\label{f:PureS}
\end{figure}

We chose a pure 1-soliton solution of the Toda lattice at time $t=0$. This ID is given by
\begin{equation}
\begin{aligned}
a_n &= 1 - \frac{1}{2}\frac{\sqrt{\left(1+e^{-2\kappa (n-1)}\right)\left(1+e^{-2\kappa (n+1)}\right)}}{1+e^{-2\kappa n}},\\
b_n &= \frac{e^{-\kappa} - e^{\kappa}}{2}\left(\frac{e^{-2\kappa n}}{1+e^{-2\kappa n}} -\frac{e^{-2\kappa (n-1)}}{1+e^{-2\kappa (n-1)}} \right),\quad \kappa=0.4,,\quad n\in\mathbb{Z}.
\end{aligned}
\label{eq:data1sol}
\end{equation}
The solution at $t = 1000$ is shown in Figure~\ref{f:PureS}.  Notice that this is the only solution we consider  without dispersive radiation.

\subsubsection{2 solitons ({\tt double})} 
\begin{figure}[h]
  \centering
  \includegraphics[width=\linewidth]{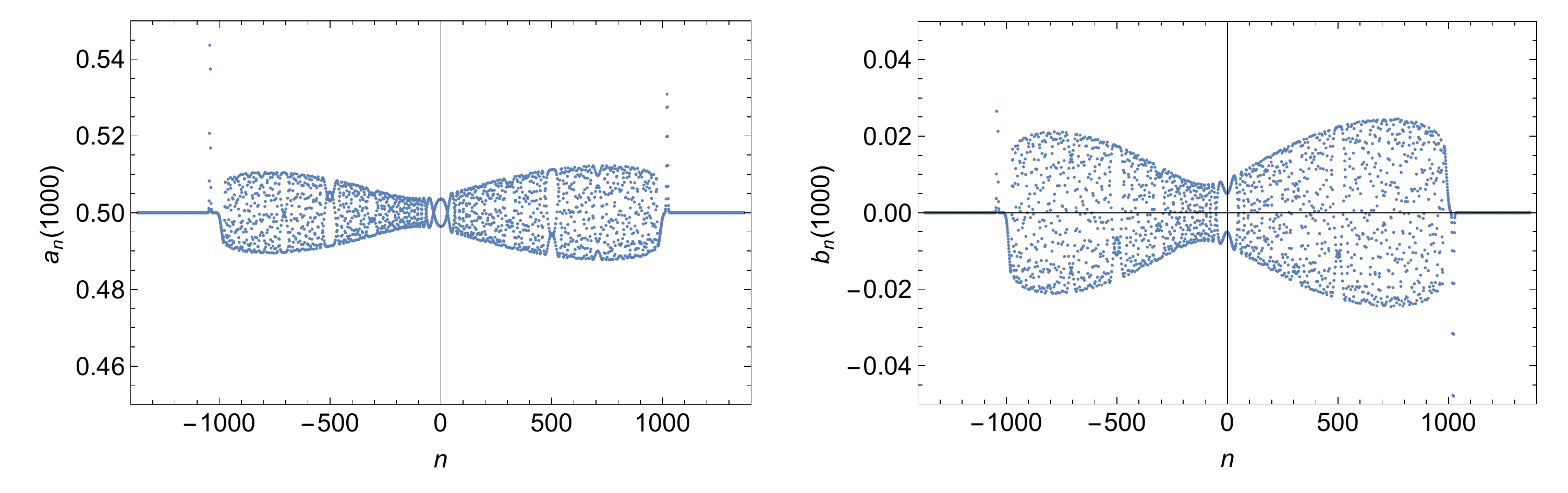}
  \caption{ The solution of the Toda lattice with {\tt double} ID at $t = 1000$.}\label{f:double}
\end{figure}

A choice of ID that gives rise to 2 solitons and radiation is
\begin{equation}
a_n = \frac{1}{2} + \frac{4}{5}n e^{-n^2}, \quad
b_n = \frac{1}{10}\sech(n),\quad n\in\mathbb{Z}.
\label{eq:data2sol}
\end{equation}
The solution at $t = 1000$ is shown in Figure~\ref{f:double}.

\subsubsection{4 solitons ({\tt quad})} 
\begin{figure}[h]
  \centering
  \includegraphics[width=\linewidth]{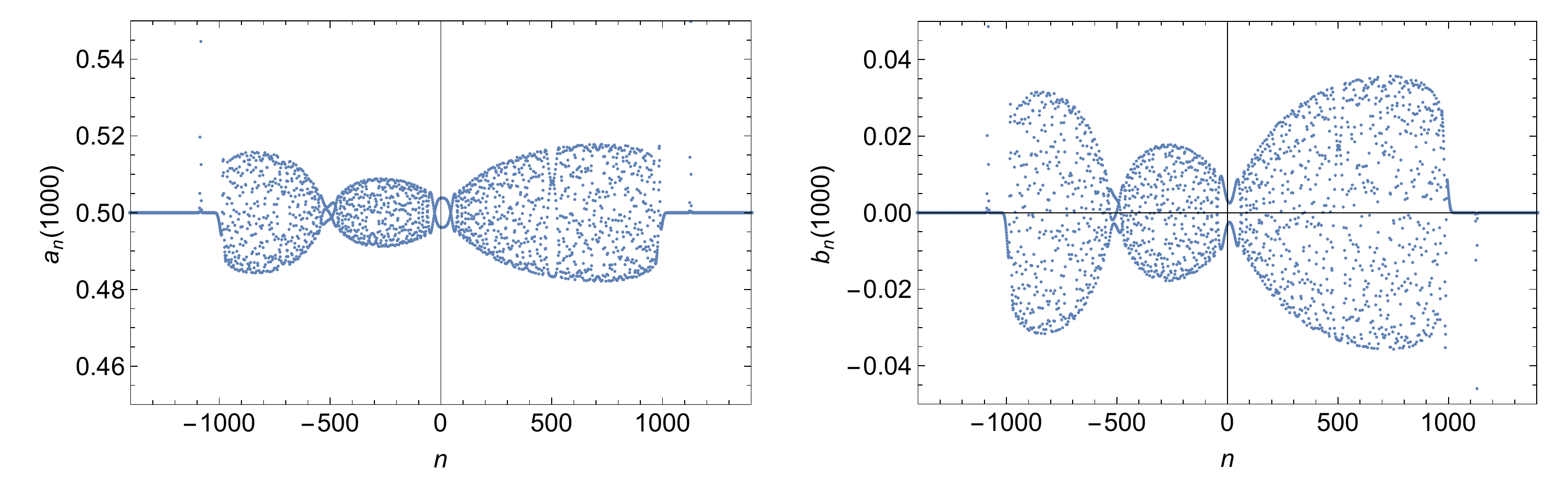}
  \caption{ The solution of the Toda lattice with {\tt quad} ID at $t = 1000$.}\label{f:quad}
\end{figure}

ID that gives rise to 4 solitons and radiation is given by
\begin{equation}
a_n =\left| \tfrac{1}{2} - n e^{-n^2 +n}\right|,\quad
b_n = n\sech(n),\quad n\in\mathbb{Z}.
\label{eq:data4sol}
\end{equation}
The solution at $t = 1000$ is shown in Figure~\ref{f:quad}.

\subsubsection{Dirac $\delta$-type ({\tt dirac})} 
\begin{figure}[h]
  \centering
  \includegraphics[width=\linewidth]{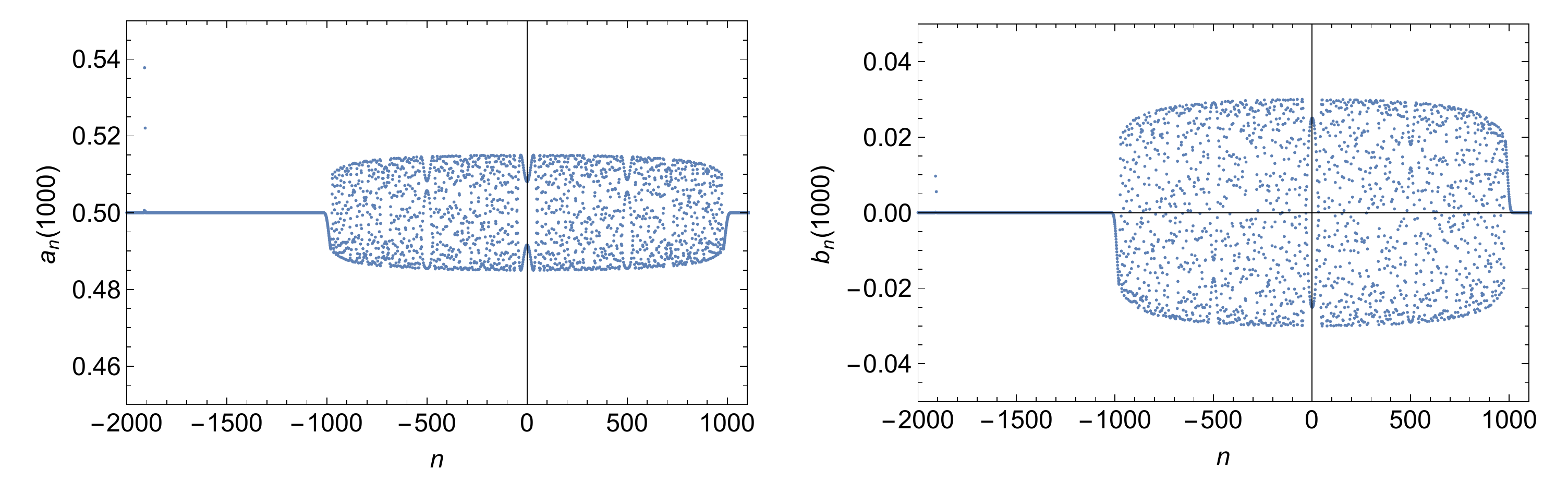}
  \caption{ The solution of the Toda lattice with {\tt dirac} ID at $t = 1000$.}\label{f:dirac}
\end{figure}

Dirac-$\delta$-type ID that leads to 1 soliton and a radiating tail that is highly oscillatory is simply given by
\begin{equation}
a_n = \frac{1}{2},\quad
b_n = \begin{cases} 4, & n= 0\\ 0, & \text{otherwise} \end{cases}\,,\quad n\in\mathbb{Z}.
\label{eq:datadirac}
\end{equation}
The solution at $t = 1000$ is shown in Figure~\ref{f:dirac}.

\begin{figure}[tbp]
  \centering
  \includegraphics[width=0.8\linewidth]{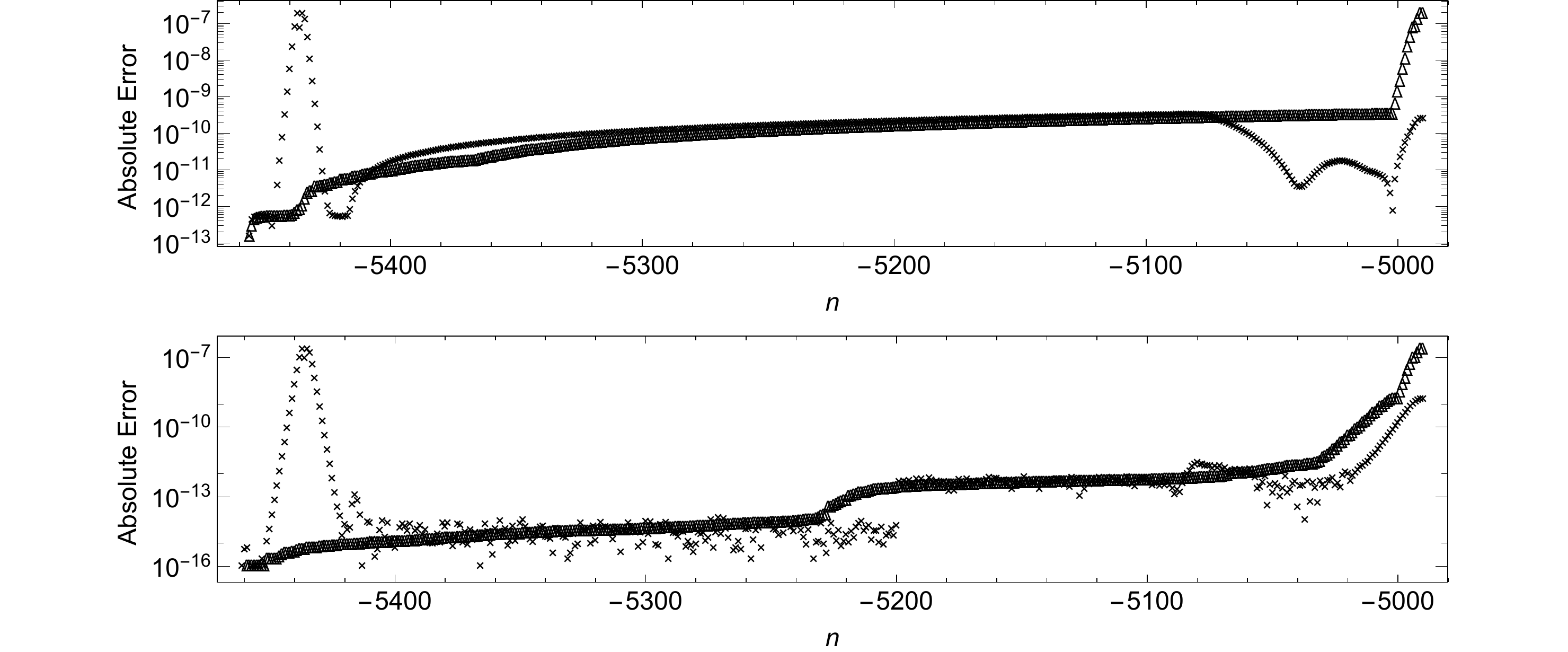}
  \caption{Absolute errors in the soliton region, for {\tt quad} ID, $dT = 0.0001$ with the {\tt ab4} (top) and {\tt sv2symp} (bottom) methods.  The $\times$ symbols represent the actual absolute errors as a function of $n$ and the $\Delta$ symbols are the errors sorted in increasing order. }
  \label{f:sortederror}
\end{figure}

\section{Comparison of Second-Order Methods}\label{s:2nd-order}

In order to compare our methods we take the following approach.  We choose three time steps $dT = 0.01, 0.001, 0.0001$, three final times $T = 1000,2000,5000$, and two regions (dispersive and soliton regions, for only $n < 0$) and examine the relative errors made in approximating the solution of the Toda lattice at the final times for every choice of ID and every choice of time-stepping method.  In accordance with the asymptotic analysis \cite{Kruger2009,Kamvissis1993,Bilman2017}, the soliton region is effectively $\mathbb Z \setminus [-T,T]$. And so, the soliton region for this work is $[-(s+100) T,-T]$ where  $s$ is the speed of the fastest moving soliton.  This can be computed from the spectrum of $\mathbf L$ \cite{Bilman2017}.  The dispersive region is $[-cT,cT]$ for $0 < c < 1$.  So, we fix our dispersive region as $[-T/2 - 50, -T/2 + 50]$.  We always take the solution computed with the numerical IST method to be our ``true'' solution.  A significant benefit of the numerical IST method is that the solution at each $(n,t)$ can be computed independently of all others.  Thus to compute the ``true'' solution we only need to compute on the intervals $[-(s+100) T,-T]$ and $[-T/2 - 50, -T/2 + 50]$ and we can refrain from computing the entire solution profile.

To define the relative error measure we use, consider the plot of errors in the soliton region, for {\tt quad} ID, $dT = 0.0001$ with the {\tt ab4} and {\tt sv2symp} methods as shown in Figure~\ref{f:sortederror}. If one uses an $\infty$-norm (or max norm) measurement (see the $\times$ symbols in Figure~\ref{f:sortederror}), one might conclude that the errors are comparable.  What this fails to account for is that the {\tt sv2symp} method has fewer errors that are near the maximum error when compared with the {\tt ab4} method.  To account for this we define the \emph{sorted norm} on $\mathbb R^n$ for $0 <d < 1$ by
\begin{align*}
  \|x\|_{\text{sort},d} = \|(y_1,y_2,\ldots,y_{\lceil d n \rceil})^T\|_2, \quad y = \mathrm{sort}(|x|).
\end{align*}
Here the function $\mathrm{sort}(\cdot)$ sorts the positive vector $|x|$ in decreasing order, and the sorted norm takes the $\ell^2$ norm of the largest $\sim dn$ entries of the vector $|x|$.  In all our computation we take $d = 0.1$, taking 10\% of the entries. This is a hybrid of the $\ell^2$-norm and the $\infty$-norm Then the error of a vector $x$ relative to a vector $y$ with background $c$ is given by
\begin{align}\label{e:rel}
  \mathrm{rel}_{y,c}(x) = \frac{\|x - y\|_{\text{sort},d}}{\|c - y\|_{\text{sort},d}}.
\end{align}
We introduce the background $c$ because for fixed $n$, $a_n(t) \to 1/2$ as $t \to \infty$ while $b_n(t) \to 0$ as $t \to \infty$.  We want to approximate $a_n(t) -1/2$ and $b_n(t)$.

\begin{remark} For {\tt PureS} ID in the dispersive region, the solution is, to machine precision, zero.  A relative error metric here does not make sense and we use absolute error $\|x - y\|_{\text{sort},s}$.
  \end{remark}

\subsection{Soliton Region}
\begin{figure}[tbp]
  \includegraphics[width=.32\linewidth]{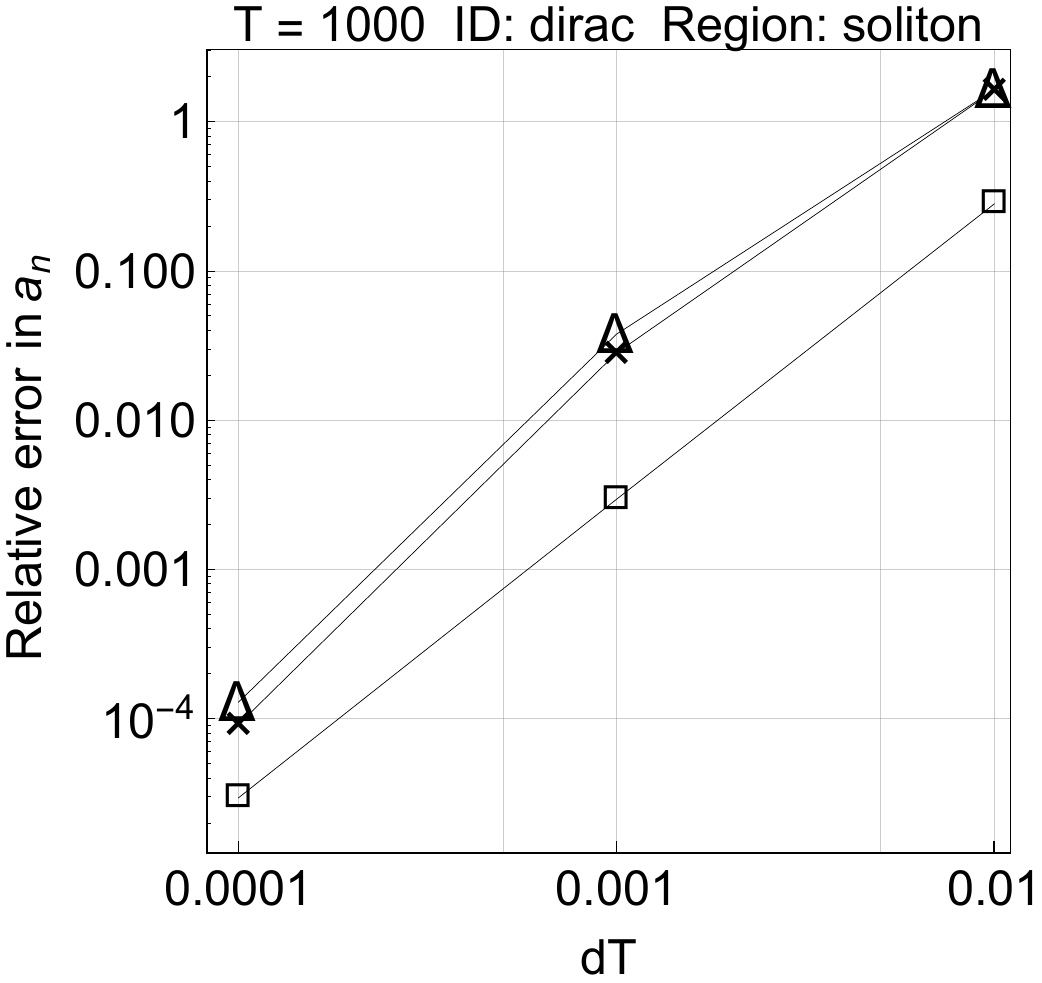}
  \includegraphics[width=.32\linewidth]{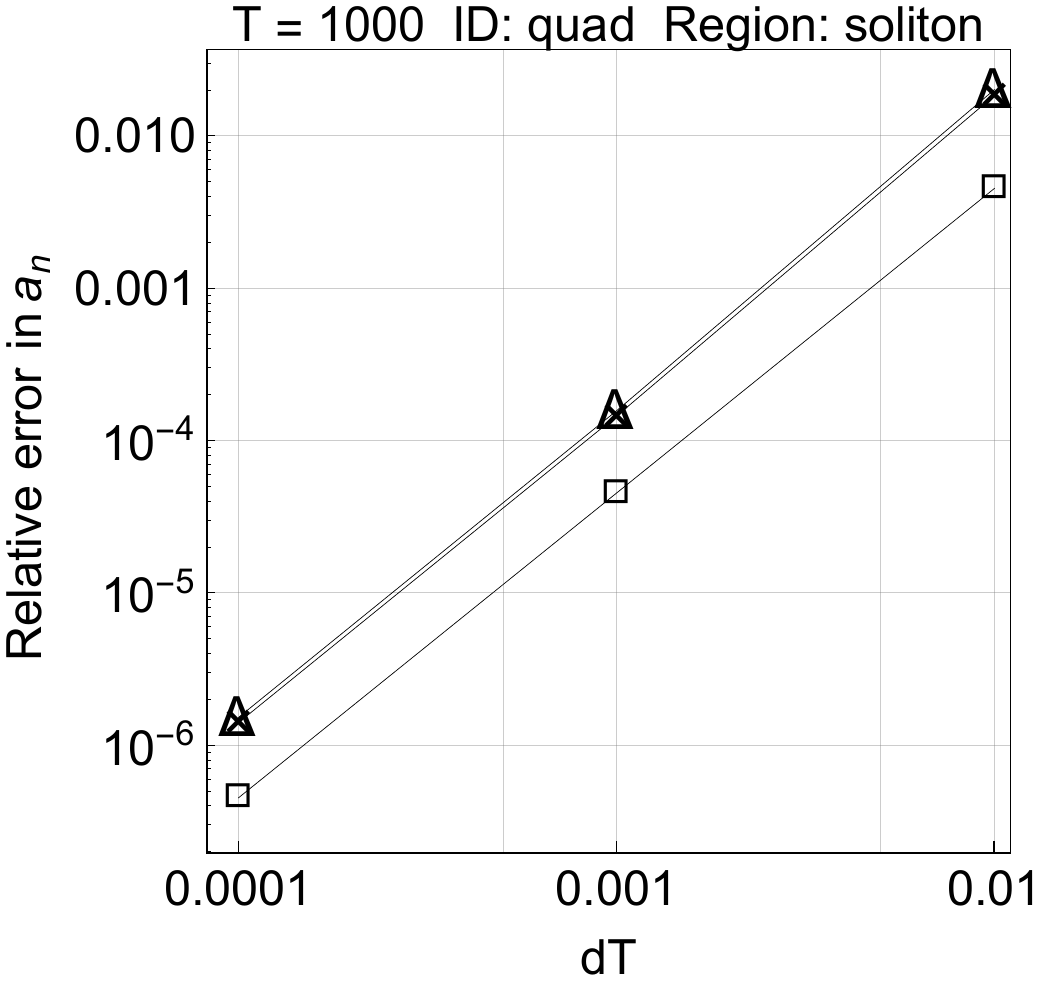}
  \includegraphics[width=.32\linewidth]{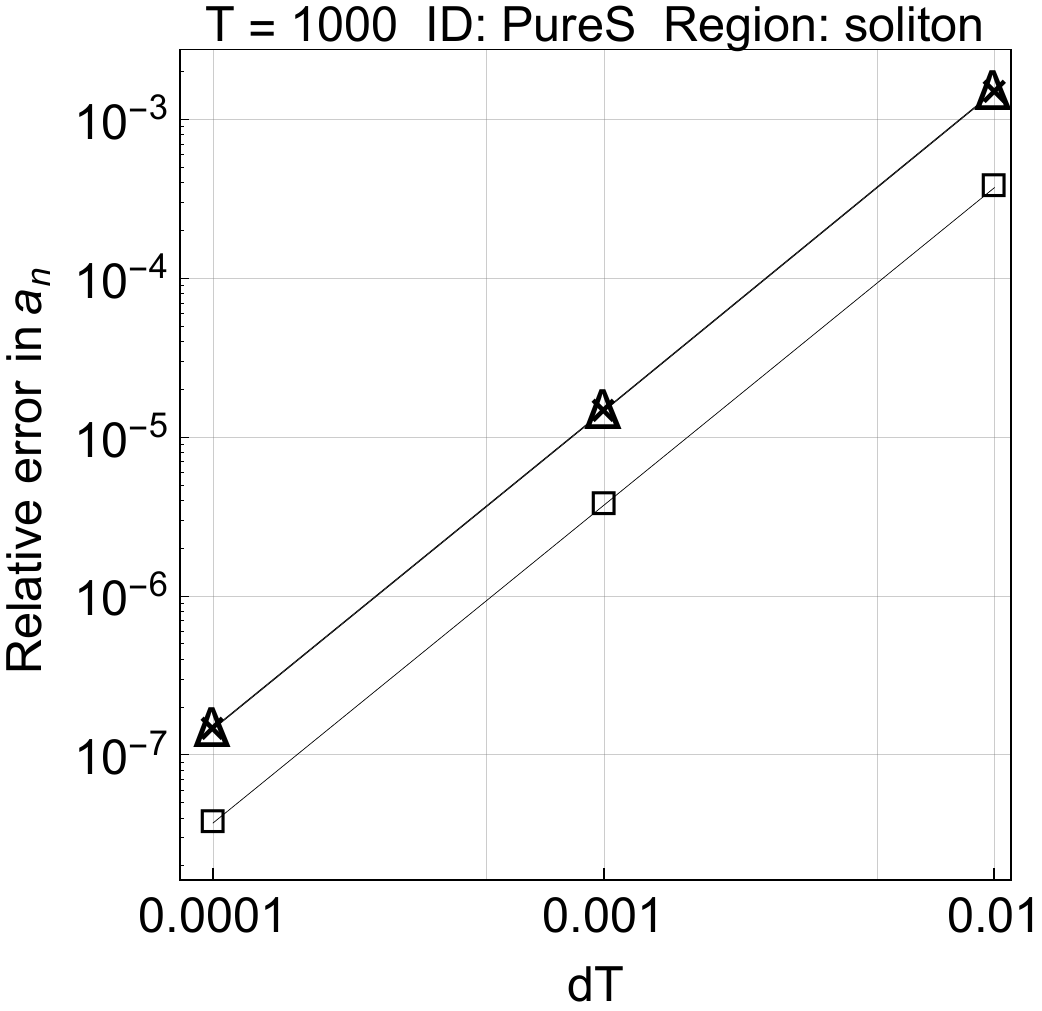}
  \caption{Relative errors in the soliton region for the second-order time-stepping methods ({\tt midpoint ($\times$), midpointqp ($\Delta$), sv2symp ($\square$)}) at $T = 1000$ plotted versus $dT$ for three choices of time step for three different choices of ID.}
    \label{f:2-soliton-1000}
\end{figure}

We first consider the errors $\mathrm{rel}_{y,1/2}(x)$ made in the approximation of the solution $a_n(T)$ of the Toda lattice in the soliton region $[-(s+100) T,-T]$ at time $T$.  Here $y$ is chosen to be the reference solution obtained by the numerical IST method.  In each panel of Figure~\ref{f:2-soliton-1000} we plot the relative error of the computed solution plotted versus $dT$.  In all panels {\tt sv2symp} out performs the other methods.

\begin{figure}[tbp]
  \includegraphics[width=.32\linewidth]{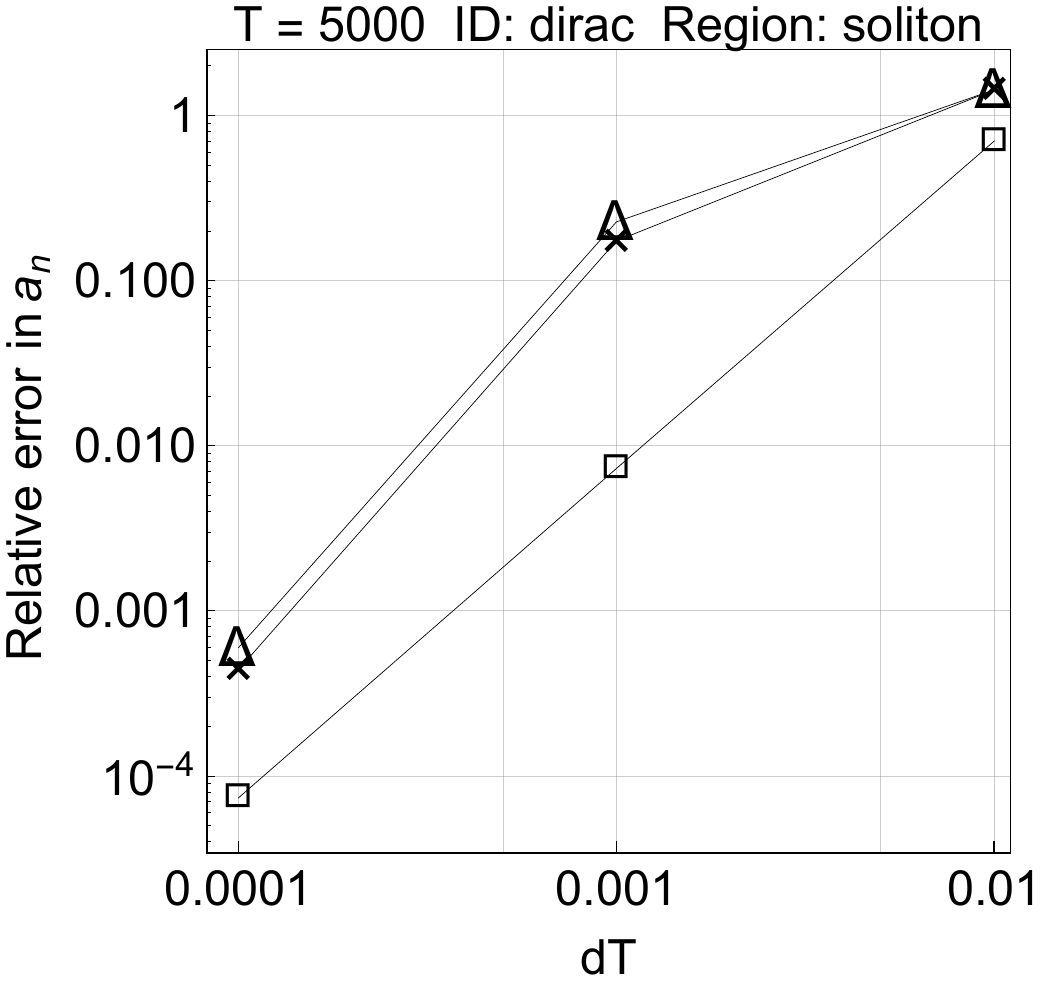}
  \includegraphics[width=.32\linewidth]{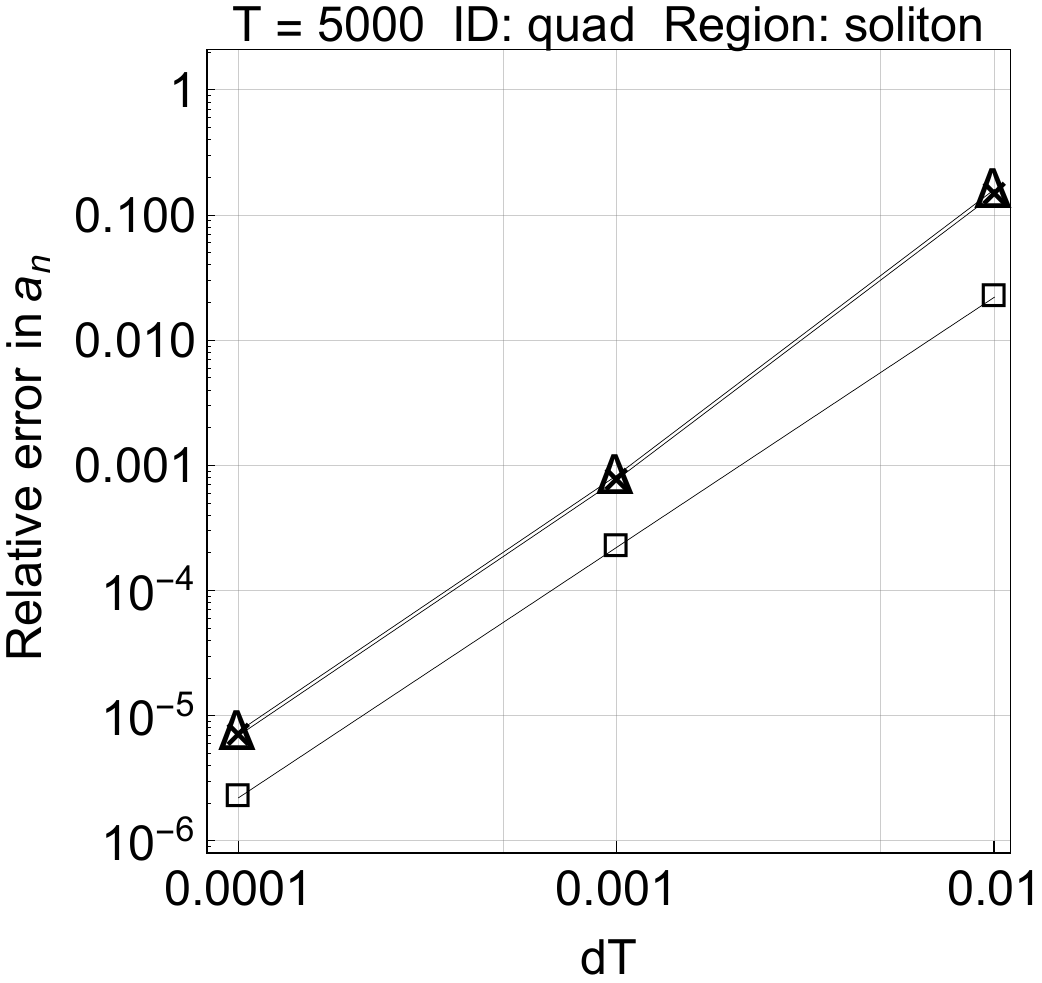}
  \includegraphics[width=.32\linewidth]{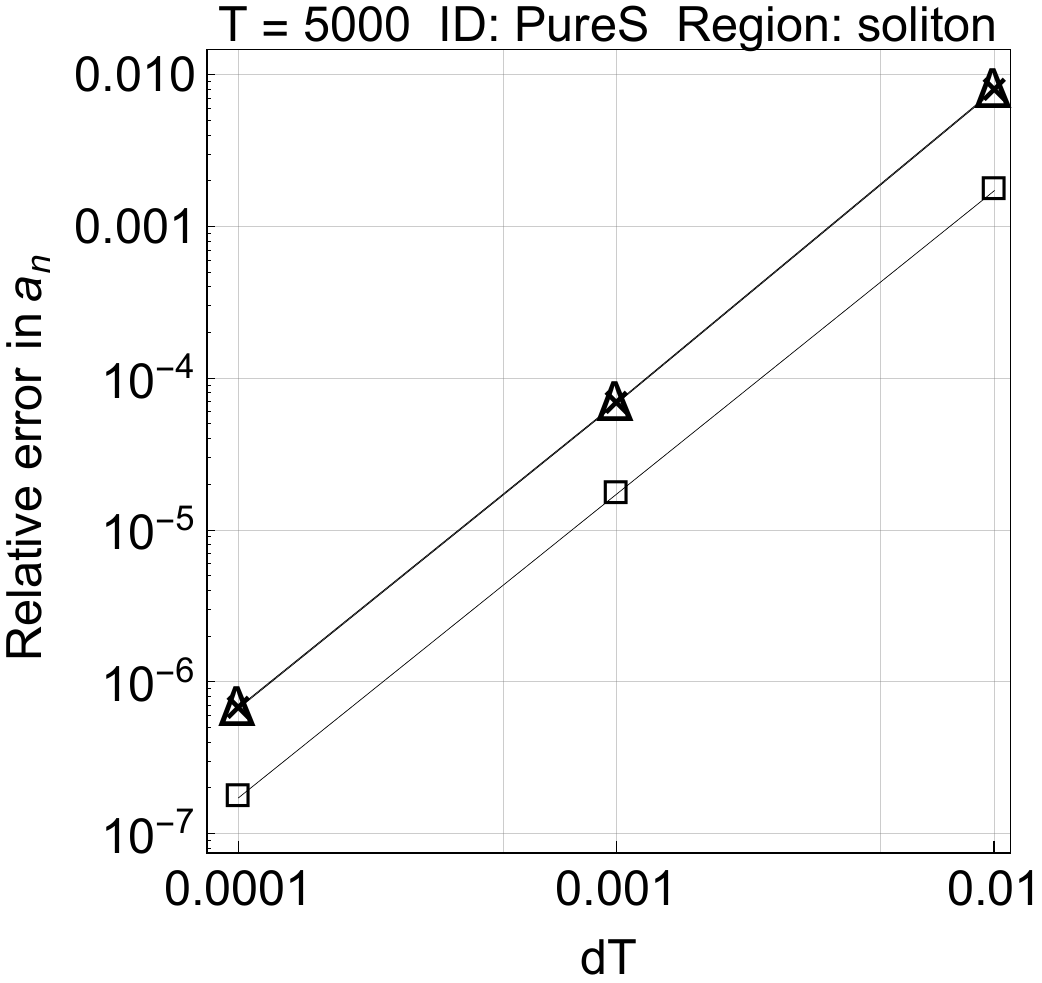}
  \caption{Relative errors in the soliton region for the second-order time-stepping methods ({\tt midpoint ($\times$), midpointqp ($\Delta$), sv2symp ($\square$)}) at $T = 5000$ plotted versus $dT$ for three choices of time step for three different choices of ID.}
    \label{f:2-soliton-5000}
\end{figure}

\subsection{Dispersive Region}

\begin{figure}[tbp]
  \includegraphics[width=.32\linewidth]{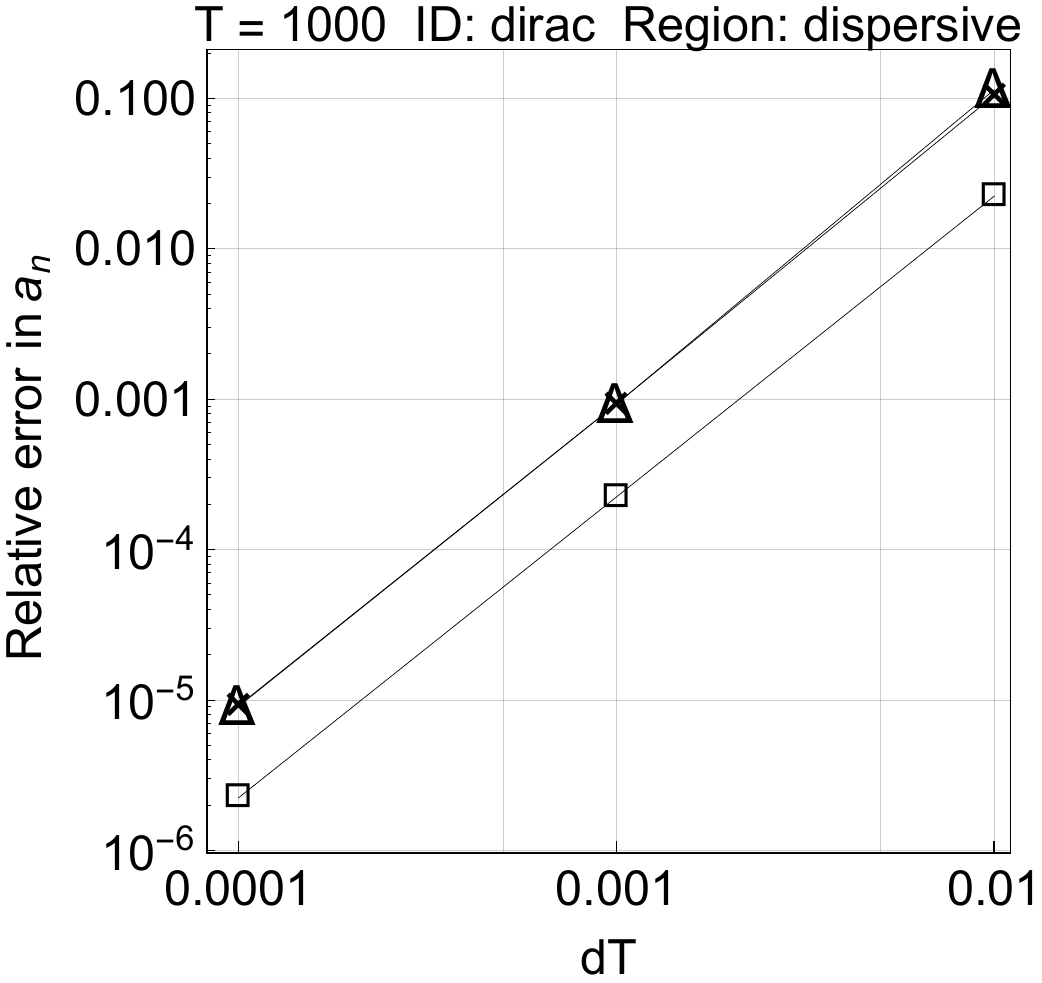}
  \includegraphics[width=.32\linewidth]{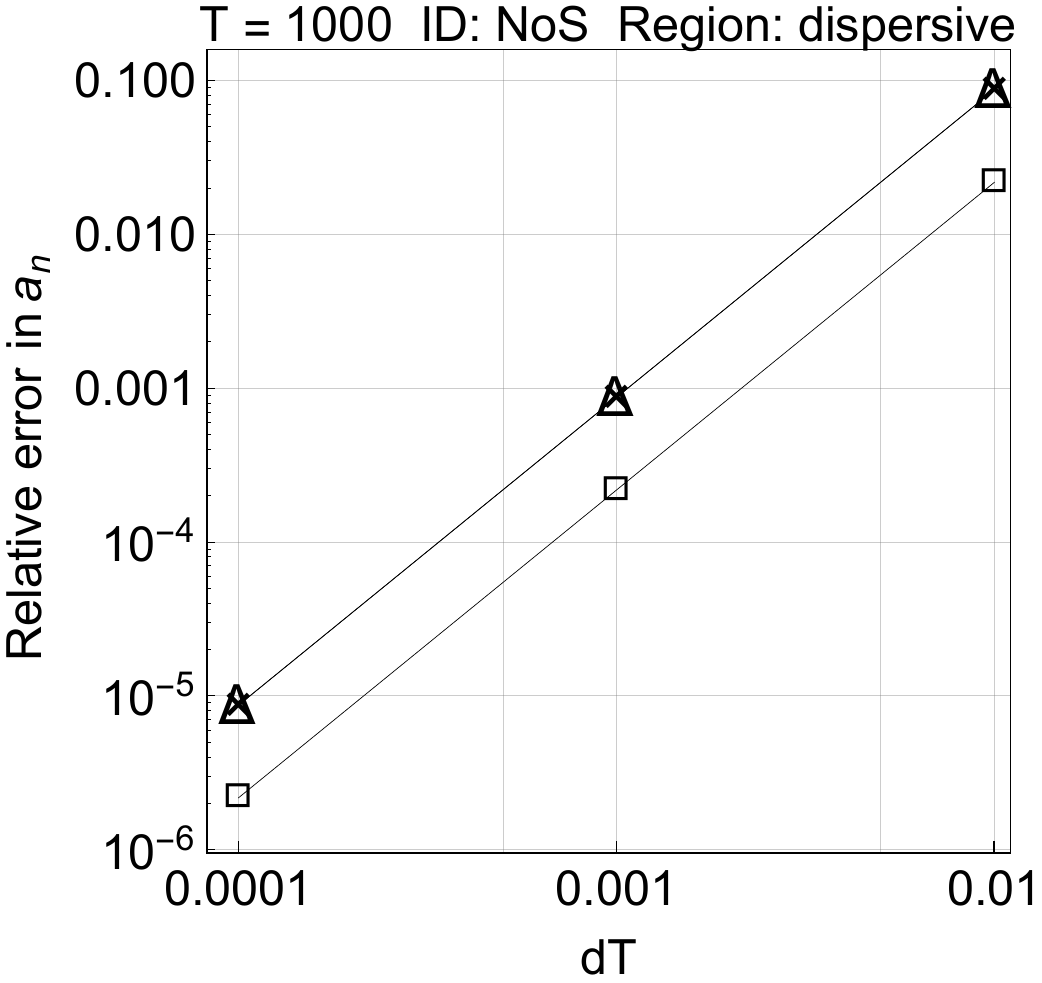}
  \includegraphics[width=.32\linewidth]{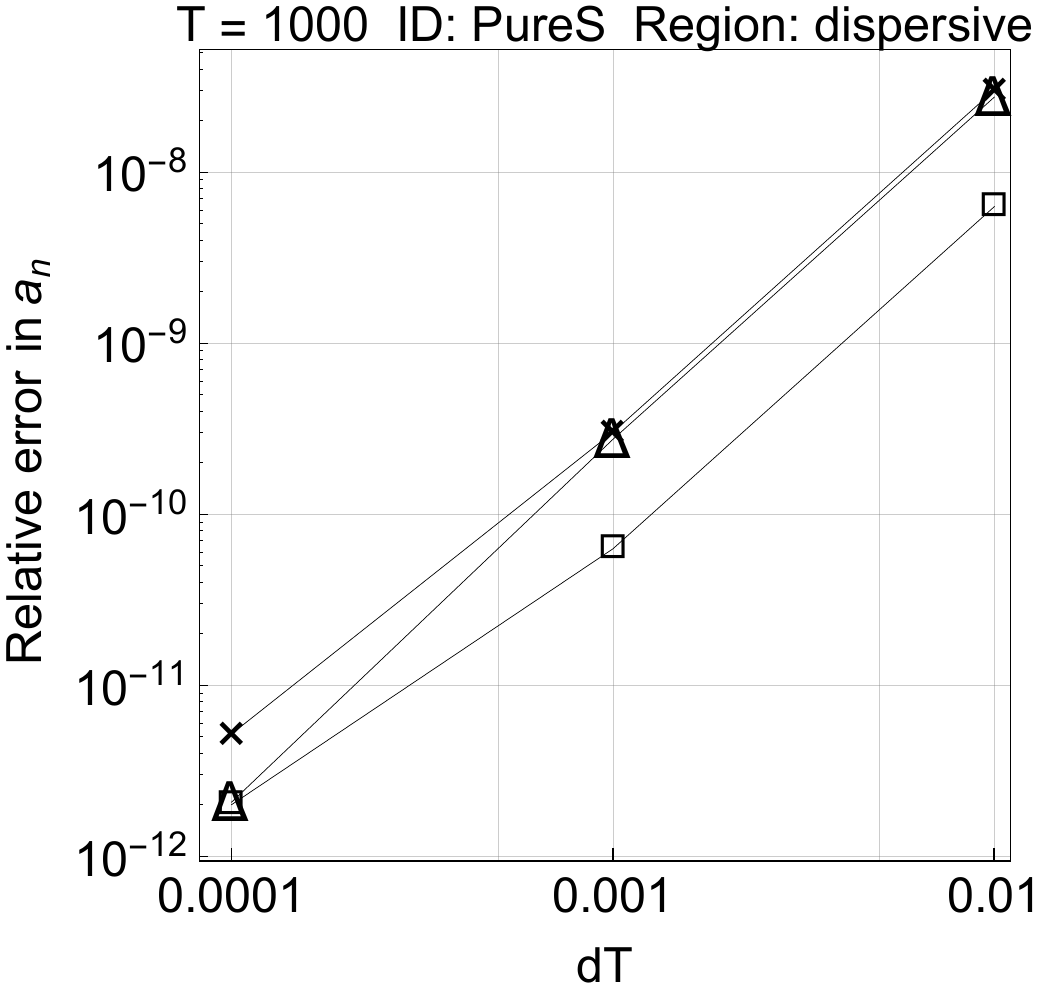}
  \caption{Relative errors in the dispersive region for the second-order time-stepping methods ({\tt midpoint ($\times$), midpointqp ($\Delta$), sv2symp ($\square$)}) at $T = 1000$ plotted versus $dT$ for three choices of time step for three different choices of ID.}
    \label{f:2-disp-1000}
\end{figure}

\begin{figure}[tbp]
  \includegraphics[width=.32\linewidth]{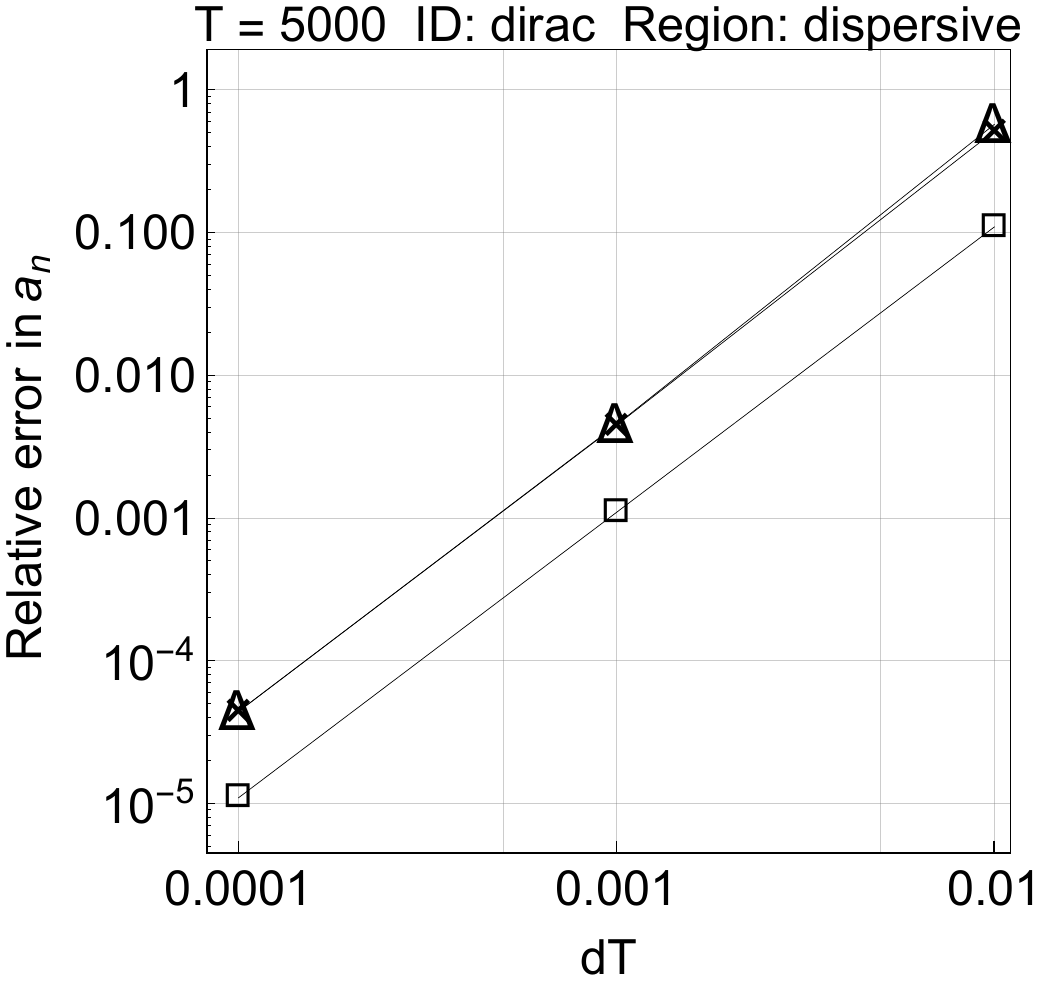}
  \includegraphics[width=.32\linewidth]{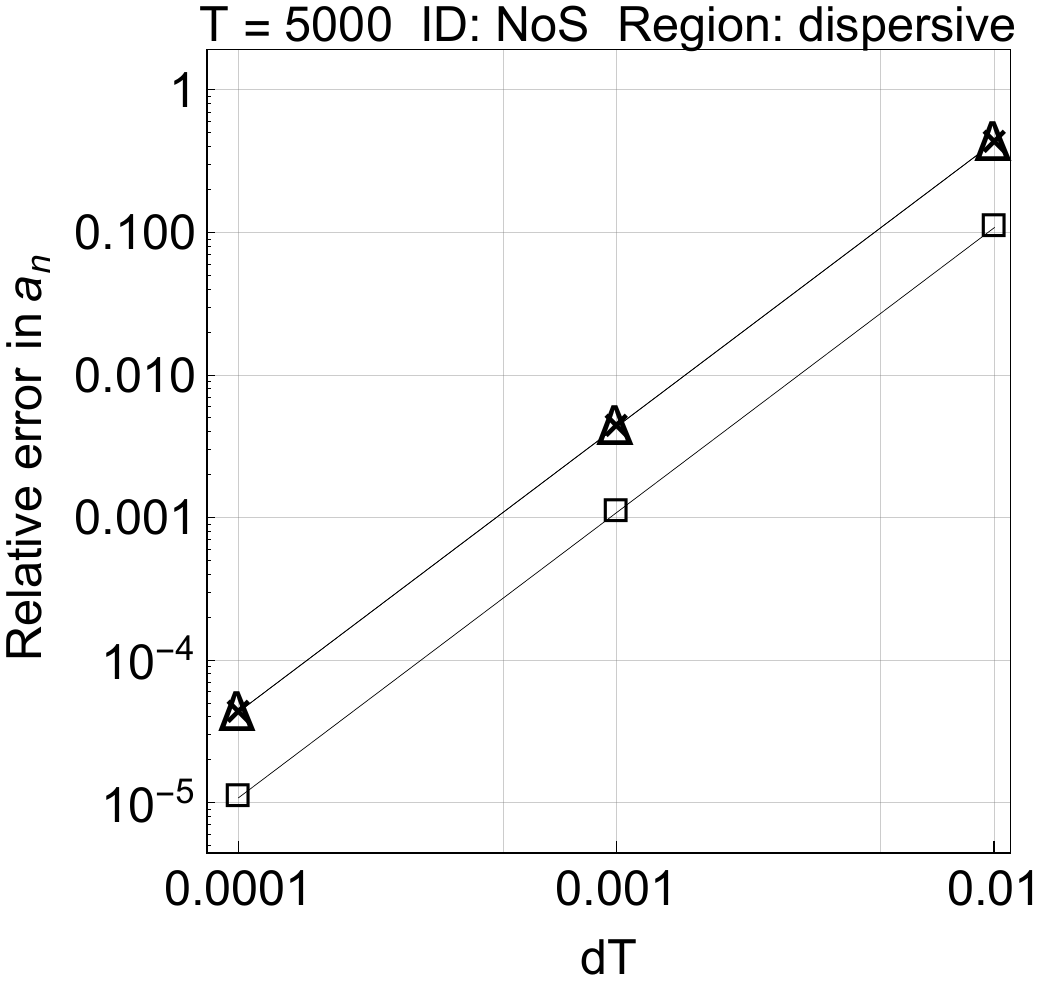}
  \includegraphics[width=.32\linewidth]{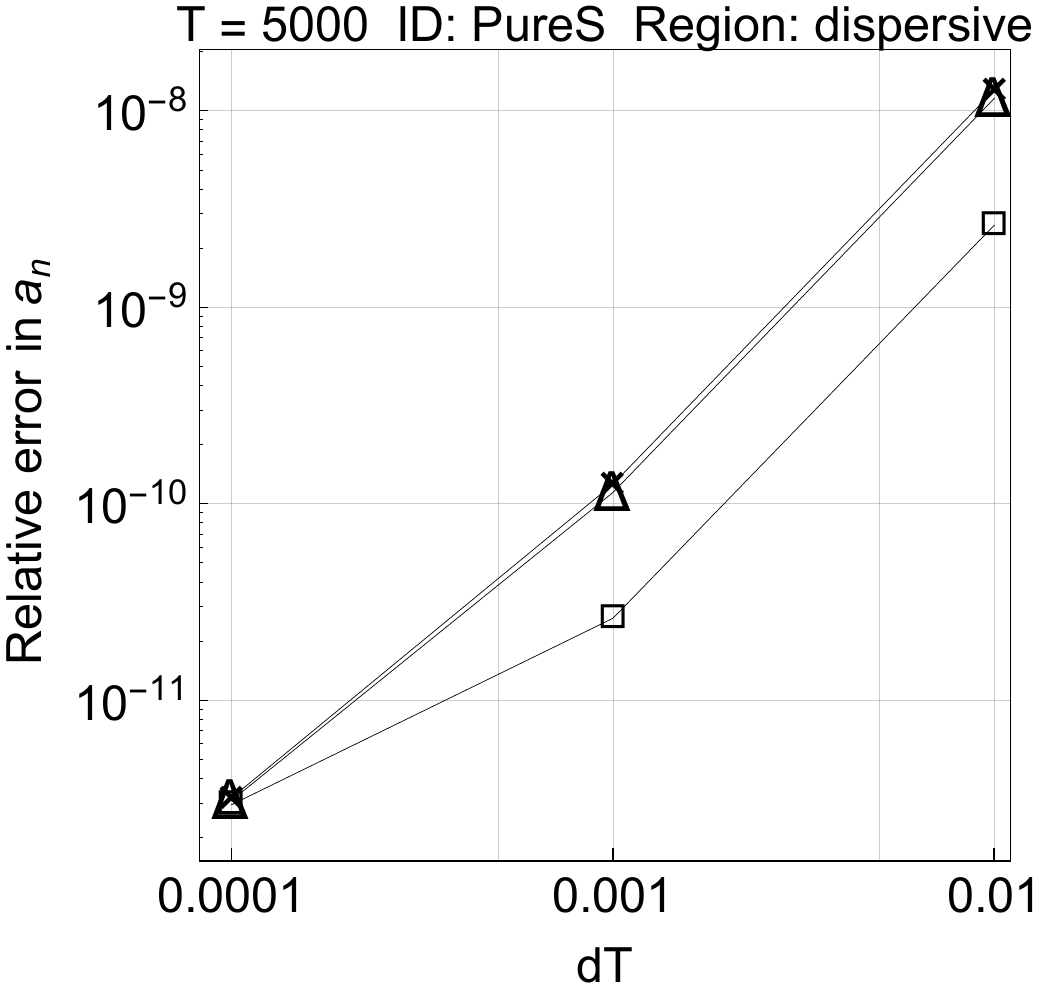}
  \caption{Relative errors in the dispersive region for the second-order time-stepping methods ({\tt midpoint ($\times$), midpointqp ($\Delta$), sv2symp ($\square$)}) at $T = 5000$ plotted versus $dT$ for three choices of time step for three different choices of ID.}
    \label{f:2-disp-5000}
\end{figure}

We now consider the errors made in the approximation of the solution $a_n(T)$ of the Toda lattice in the dispersive region region $[-T/2 - 50 ,-T/2 + 50]$ at time $T$. We can make an important point with Figure~\ref{f:2-disp-1000} and \ref{f:2-disp-5000}.  If one could only work with the pure soliton solution (right panel in the figures), that person might conclude that the {\tt midpointqp} for small enough time step performs as well as {\tt sv2symp} away from the soliton.  This is true for the {\tt PureS} ID, but not for the other ID.  This illustrates why having accurate solutions with dispersive tails to compare against is important.

\section{Comparison of Fourth-Order Methods}\label{s:higher-order}

We now move to the comparison of the fourth-order methods listed above.  We use the same relative error metric $\mathrm{rel}_{y,1/2}(x)$ as described in \eqref{e:rel}, where the reference solution $y$ is computed with the numerical IST method.

\subsection{Soliton Region}

\begin{figure}[tbp]
  \includegraphics[width=.32\linewidth]{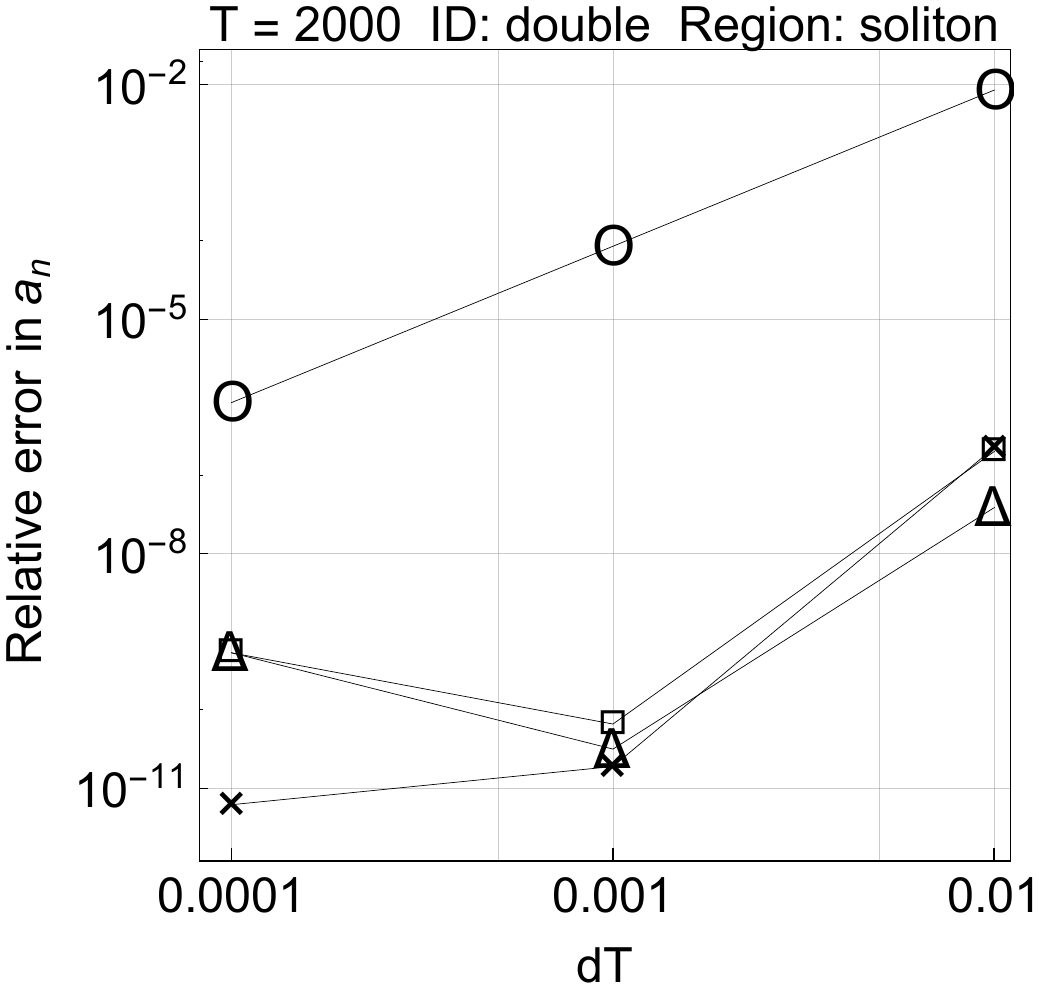}
  \includegraphics[width=.32\linewidth]{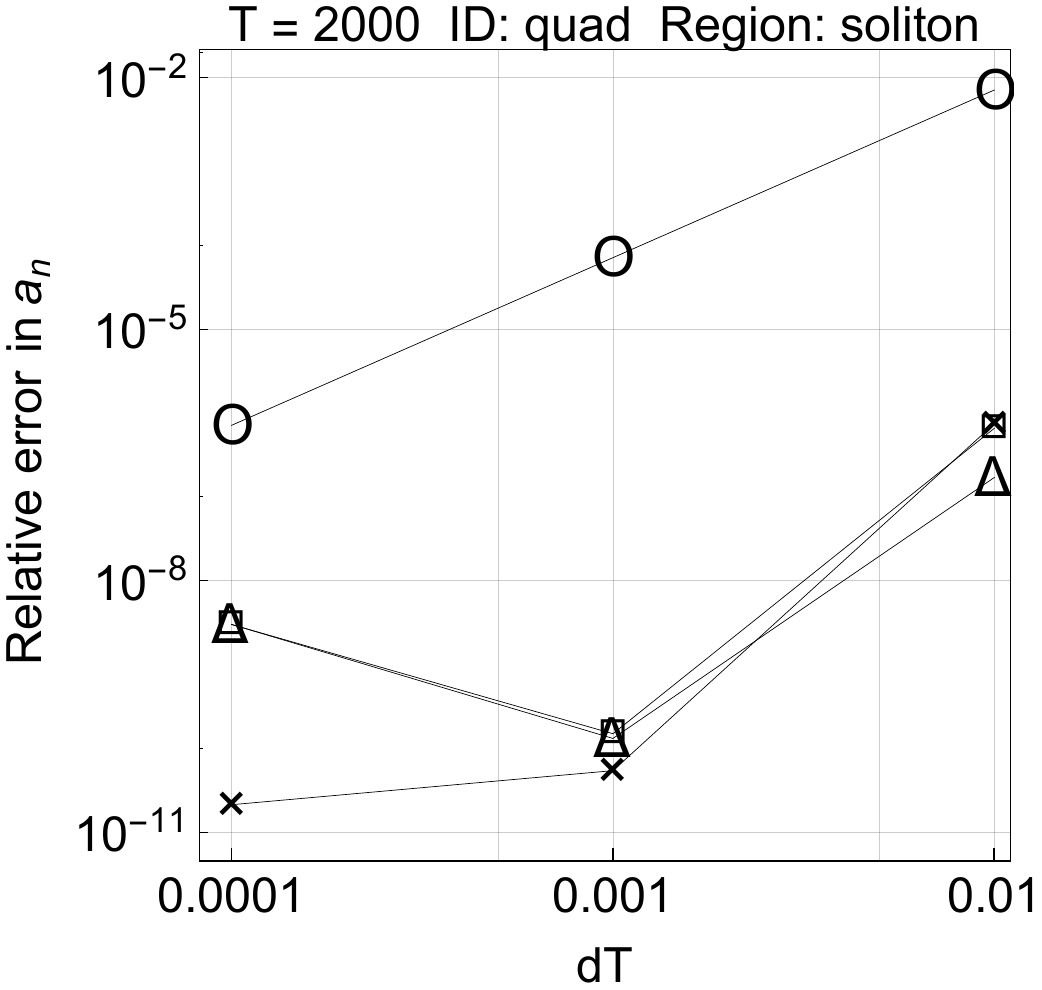}
  \includegraphics[width=.32\linewidth]{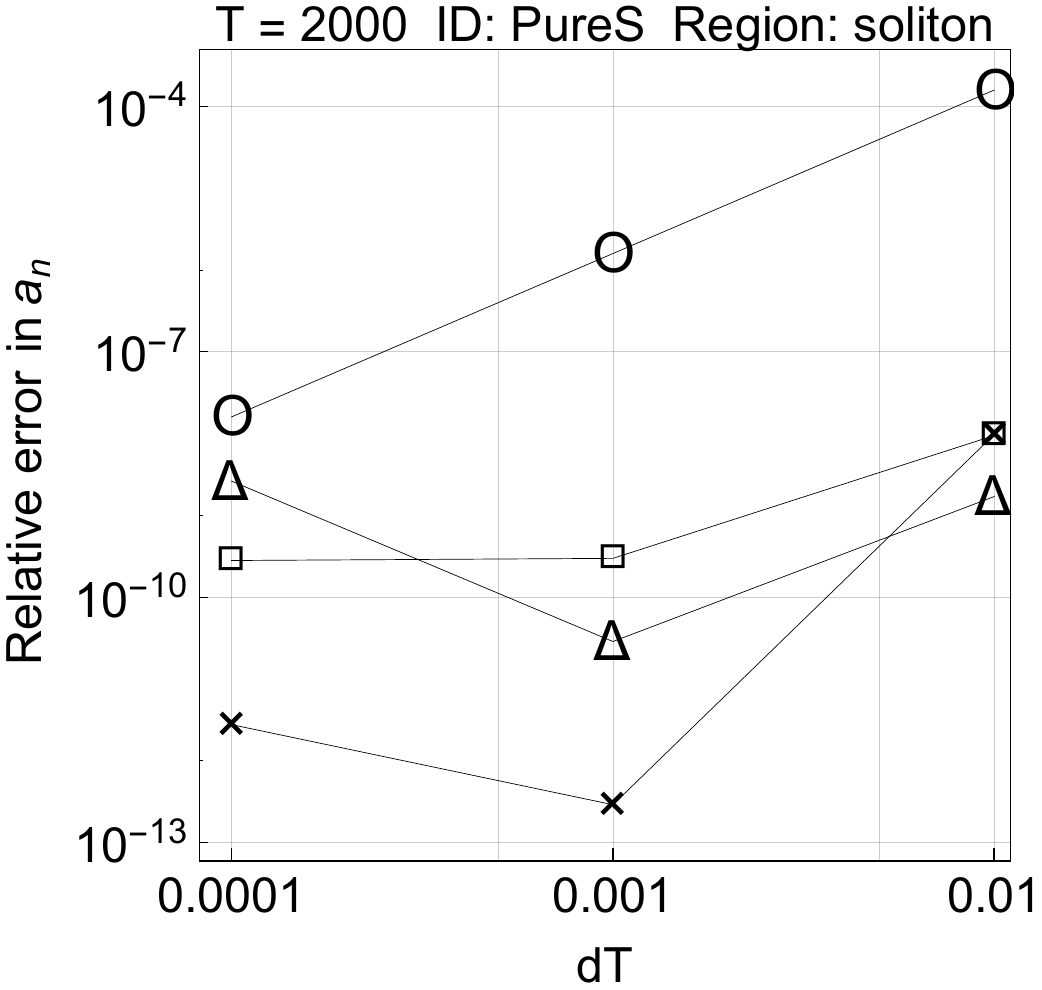}
  \caption{Relative errors in the soliton region for the fourth-order time-stepping methods ({\tt rk4 ($\square$), rk4qp ($\times$), rkf45 ($\Delta$), ab4 ($\ocircle$)}) at $T = 2000$ plotted versus $dT$ for three choices of time step for three different choices of ID.}
    \label{f:4-soliton-2000}
\end{figure}

\begin{figure}[tbp]
  \includegraphics[width=.32\linewidth]{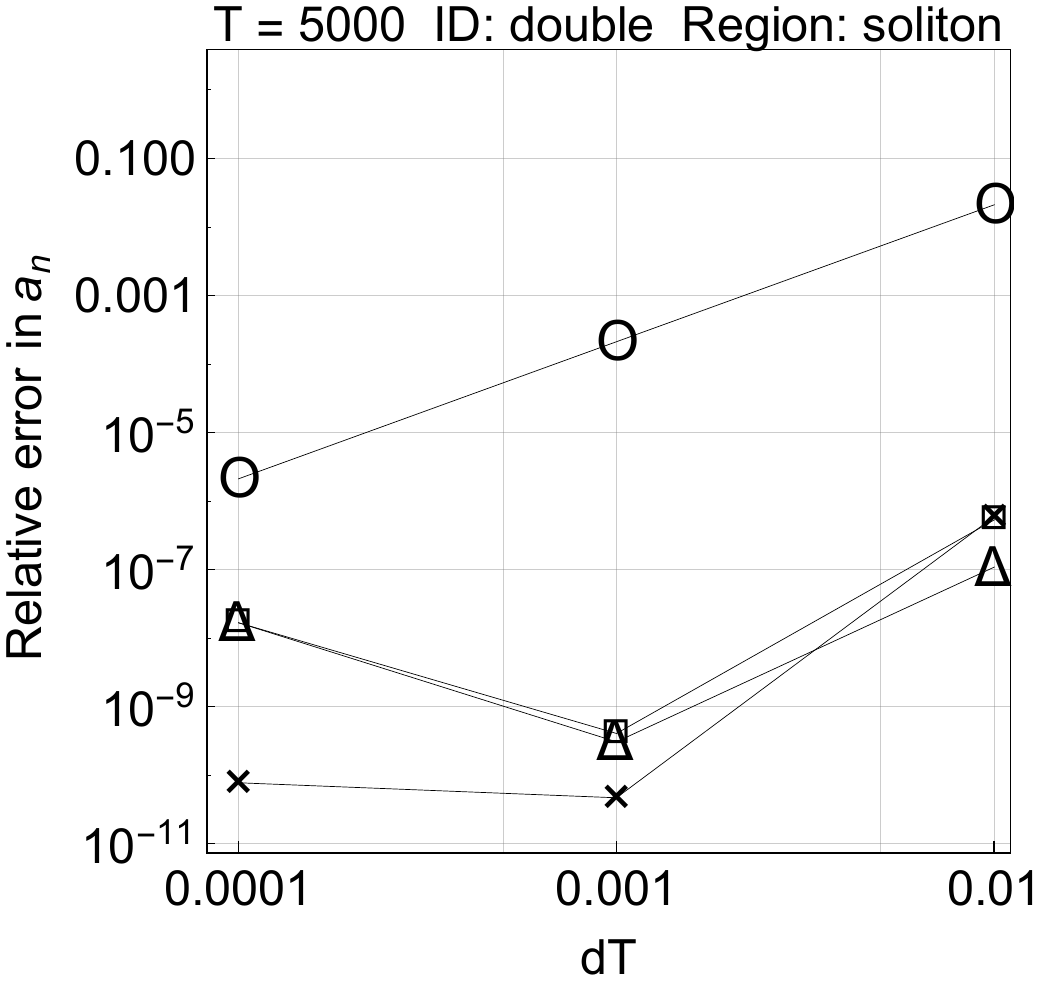}
  \includegraphics[width=.32\linewidth]{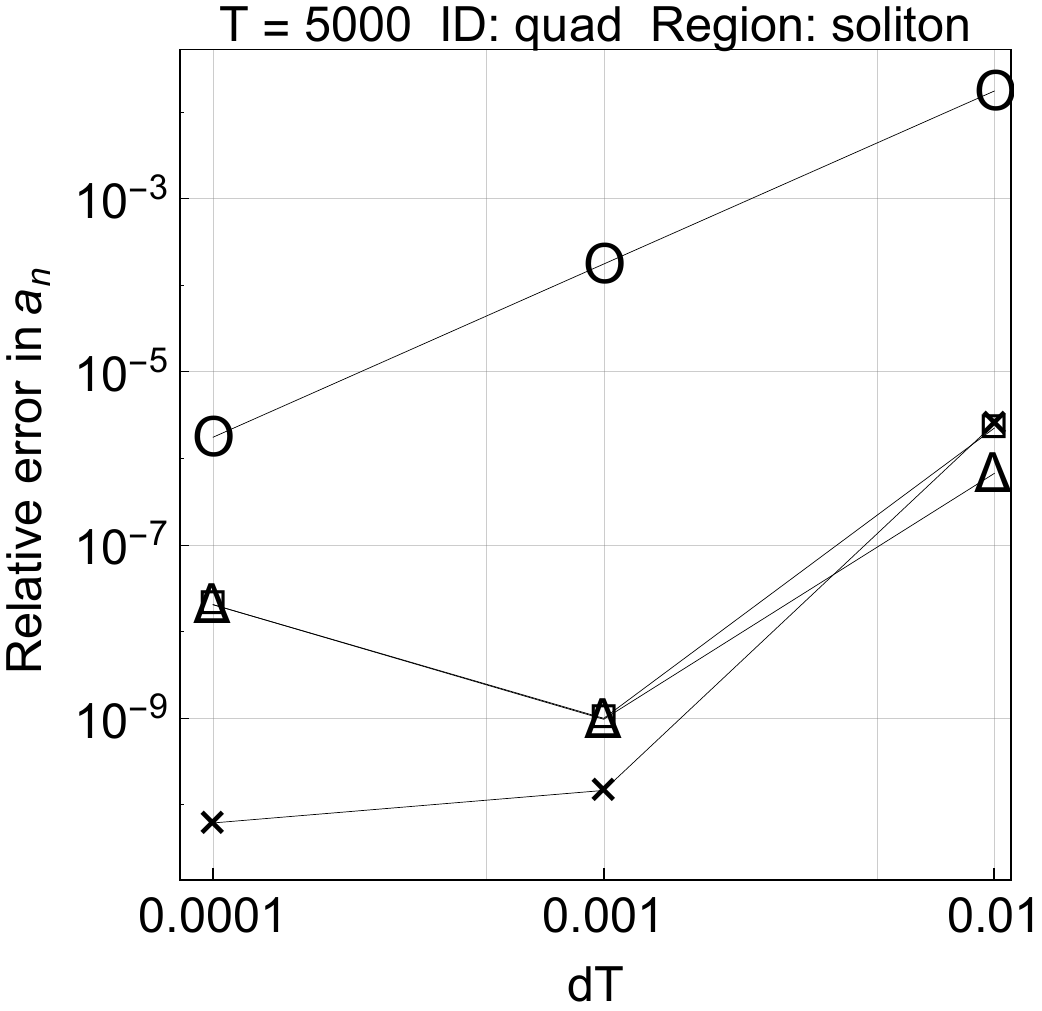}
  \includegraphics[width=.32\linewidth]{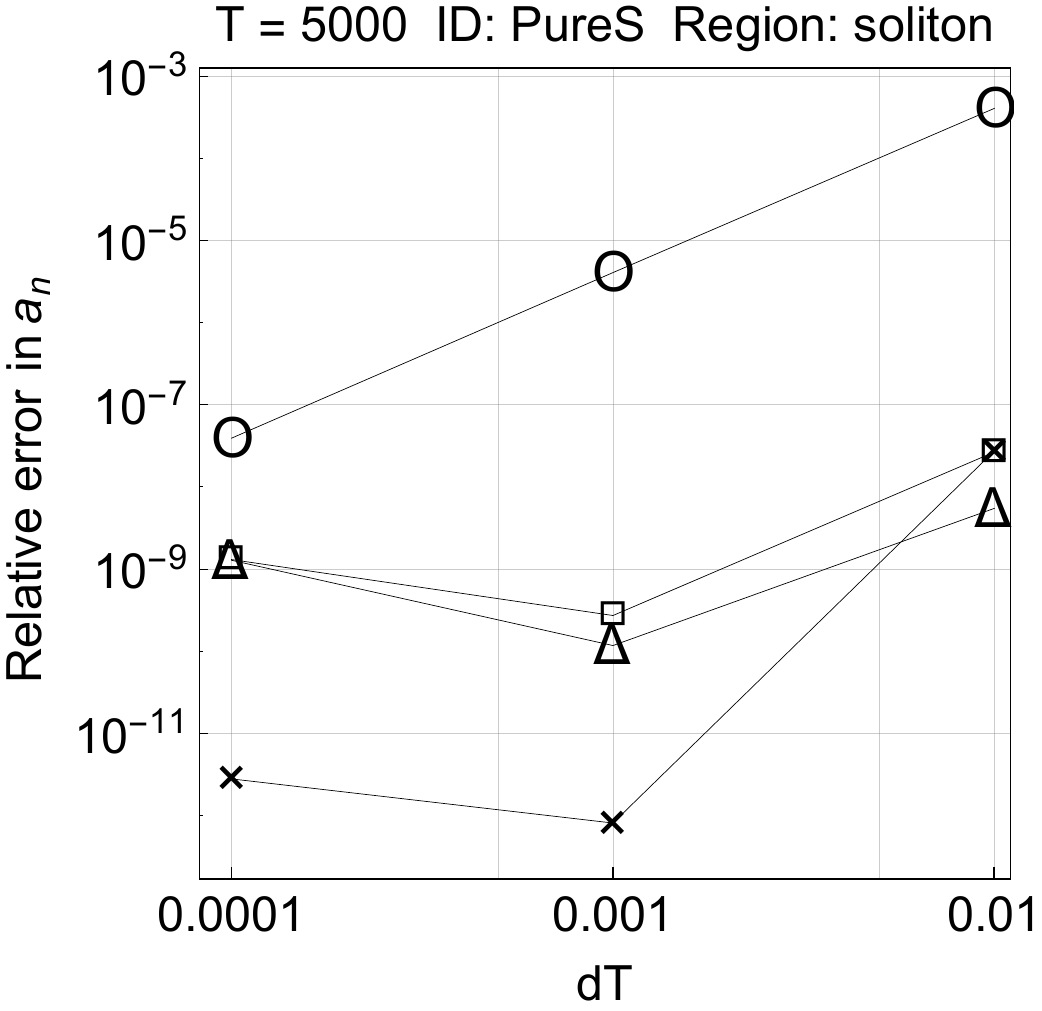}
  \caption{Relative errors in the soliton region for the fourth-order time-stepping methods ({\tt rk4 ($\square$), rk4qp ($\times$), rkf45 ($\Delta$), ab4 ($\ocircle$)}) at $T = 5000$ plotted versus $dT$ for three choices of time step for three different choices of ID.}
    \label{f:4-soliton-5000}
\end{figure}

In each of the panels of Figures~\ref{f:4-soliton-2000} and \ref{f:4-soliton-5000} we plot the relative error of the computed approximation of $a_n(T)$ plotted versus $dT$ for fourth-order methods in the soliton region.  We can see that we are operating near the maximum accuracy of these methods as the relative error can increase as $dT$ decreases.  We see that {\tt ab4} under-performs and {\tt rk4qp} is almost always the method of choice in this region.  But it is important to note that the relative error encountered for the {\tt PureS} ID is less than that encountered for the other ID.  This points again to the need for a wide class of test solutions.

\subsection{Dispersive Region}

\begin{figure}[tbp]
  \includegraphics[width=.32\linewidth]{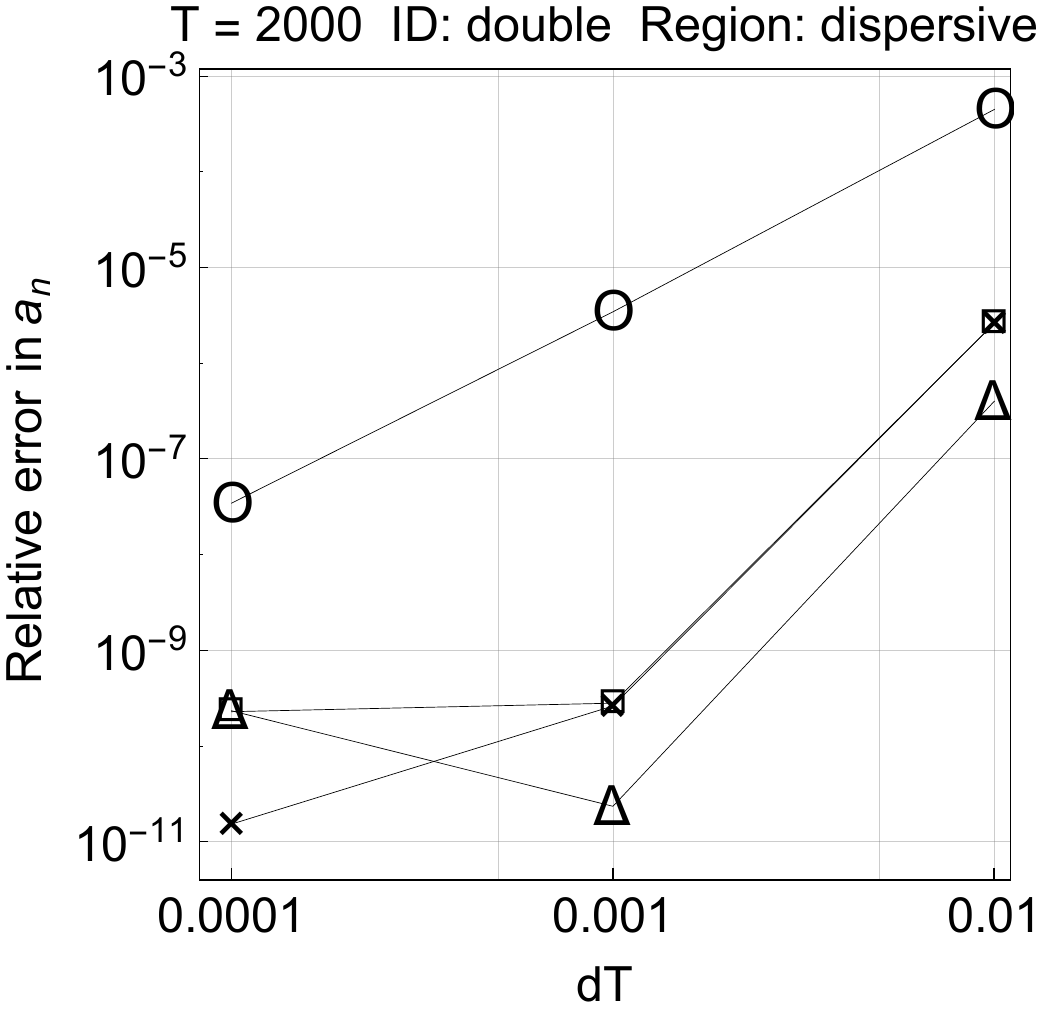}
  \includegraphics[width=.32\linewidth]{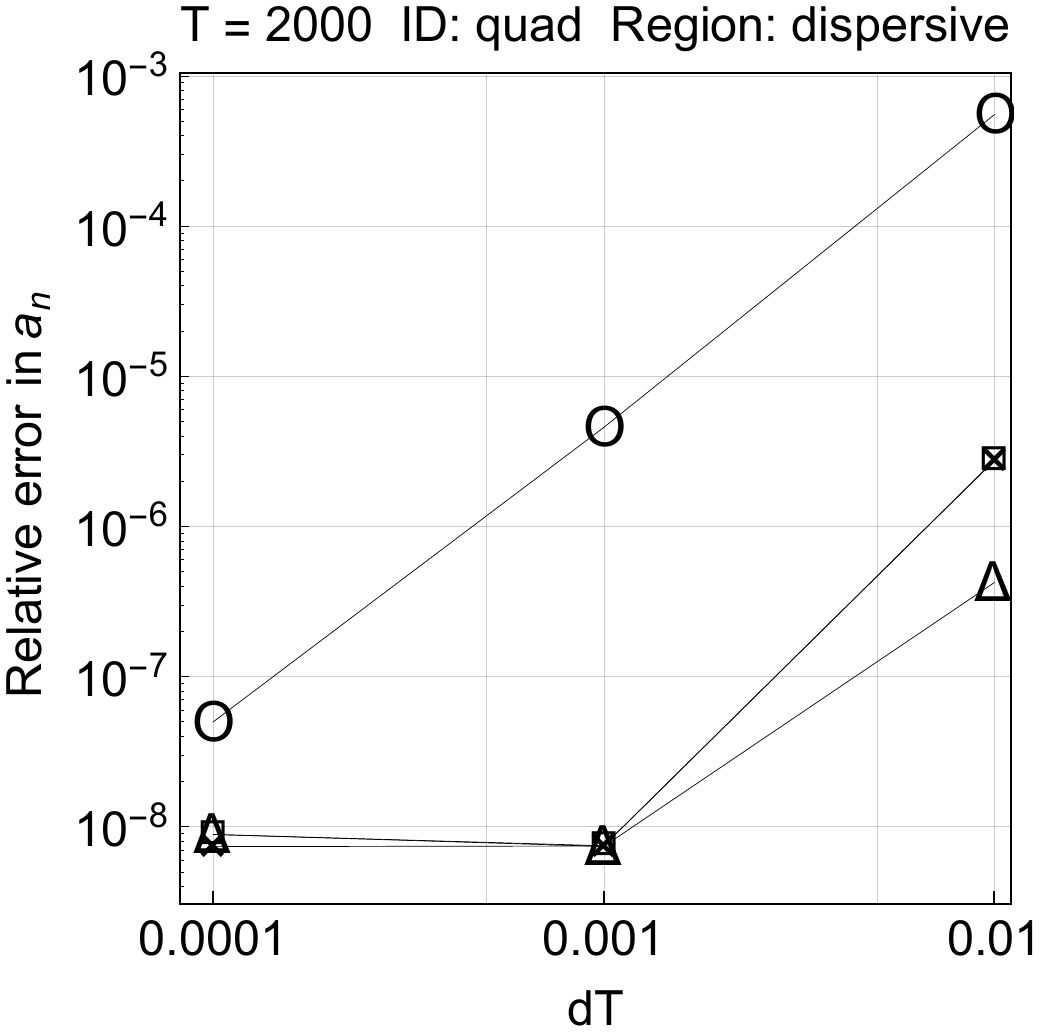}
  \includegraphics[width=.32\linewidth]{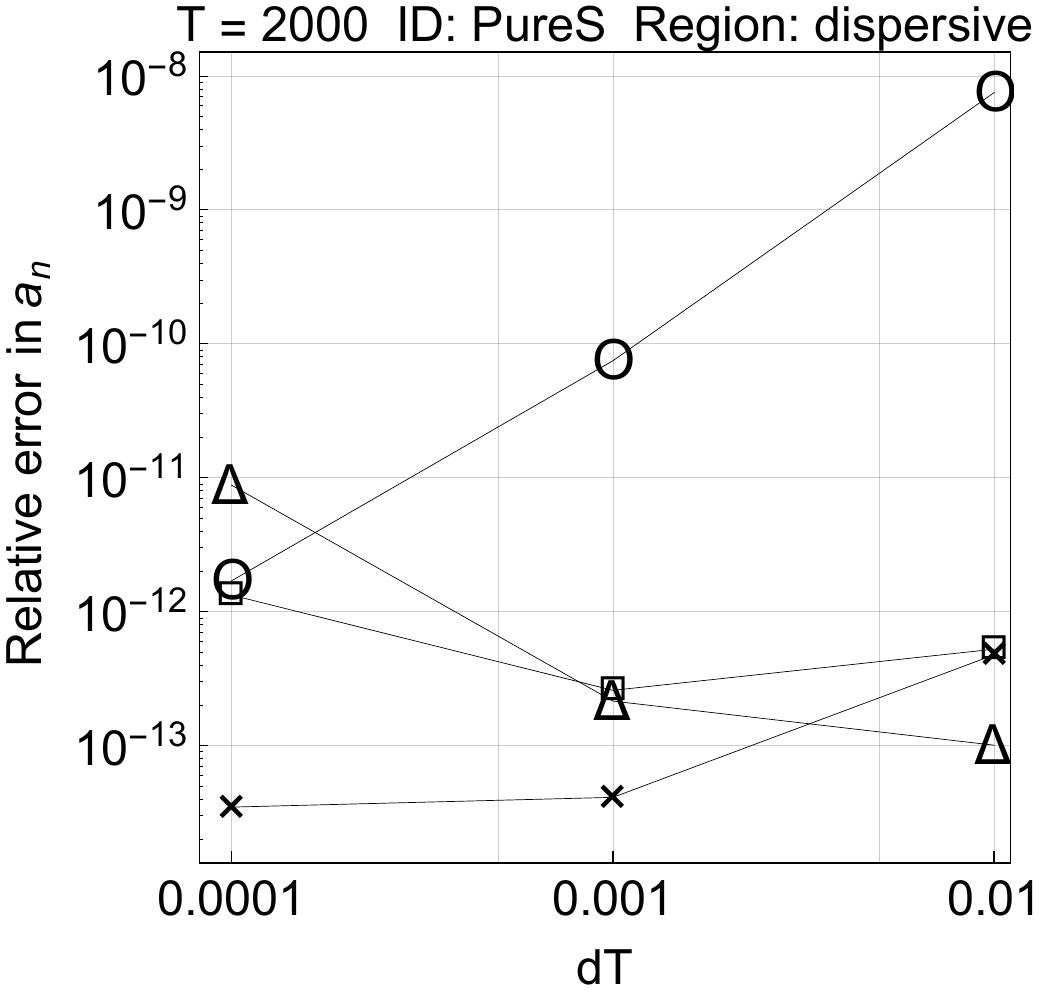}
  \caption{Relative errors in the dispersive region for the fourth-order time-stepping methods ({\tt rk4 ($\square$), rk4qp ($\times$), rkf45 ($\Delta$), ab4 ($\ocircle$)}) at $T = 2000$ plotted versus $dT$ for three choices of time step for three different choices of ID.}
    \label{f:4-disp-2000}
\end{figure}

\begin{figure}[tbp]
  \includegraphics[width=.32\linewidth]{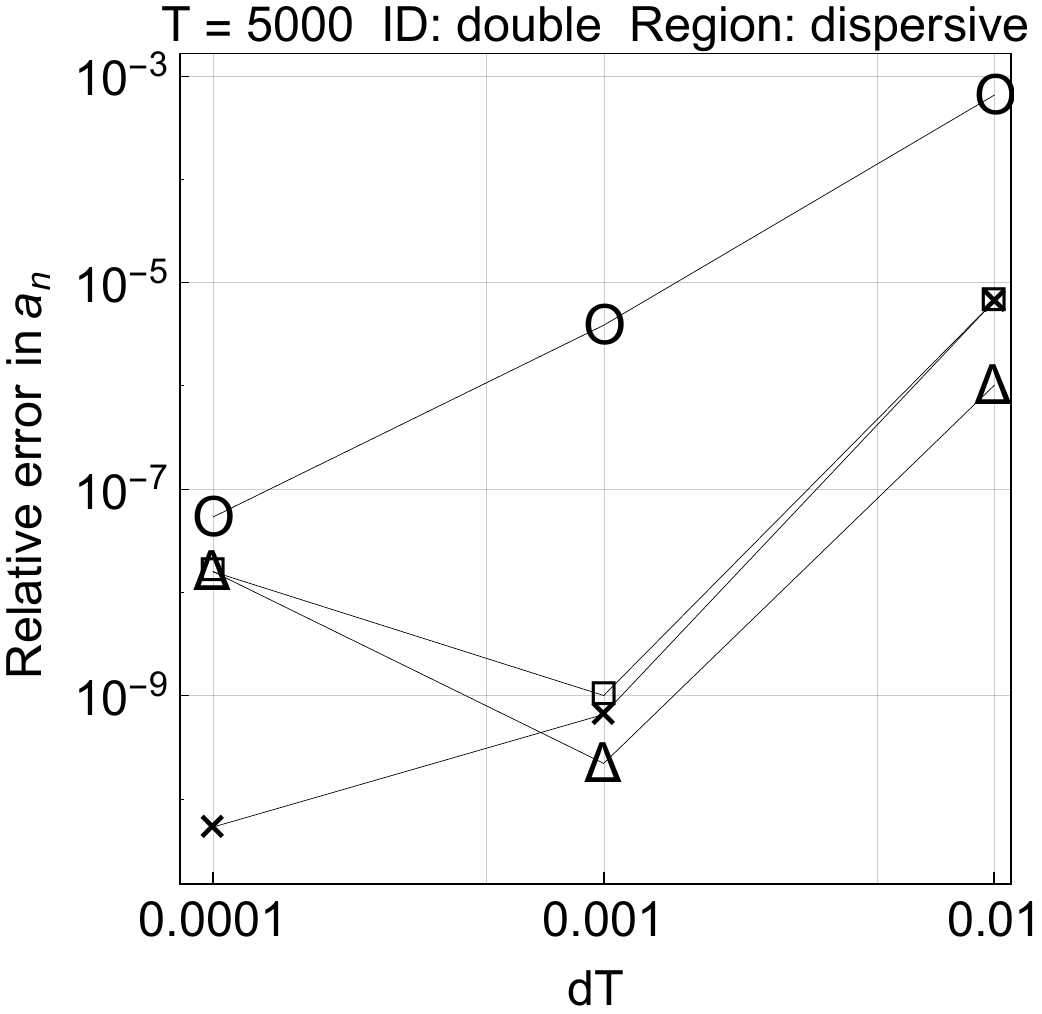}
  \includegraphics[width=.32\linewidth]{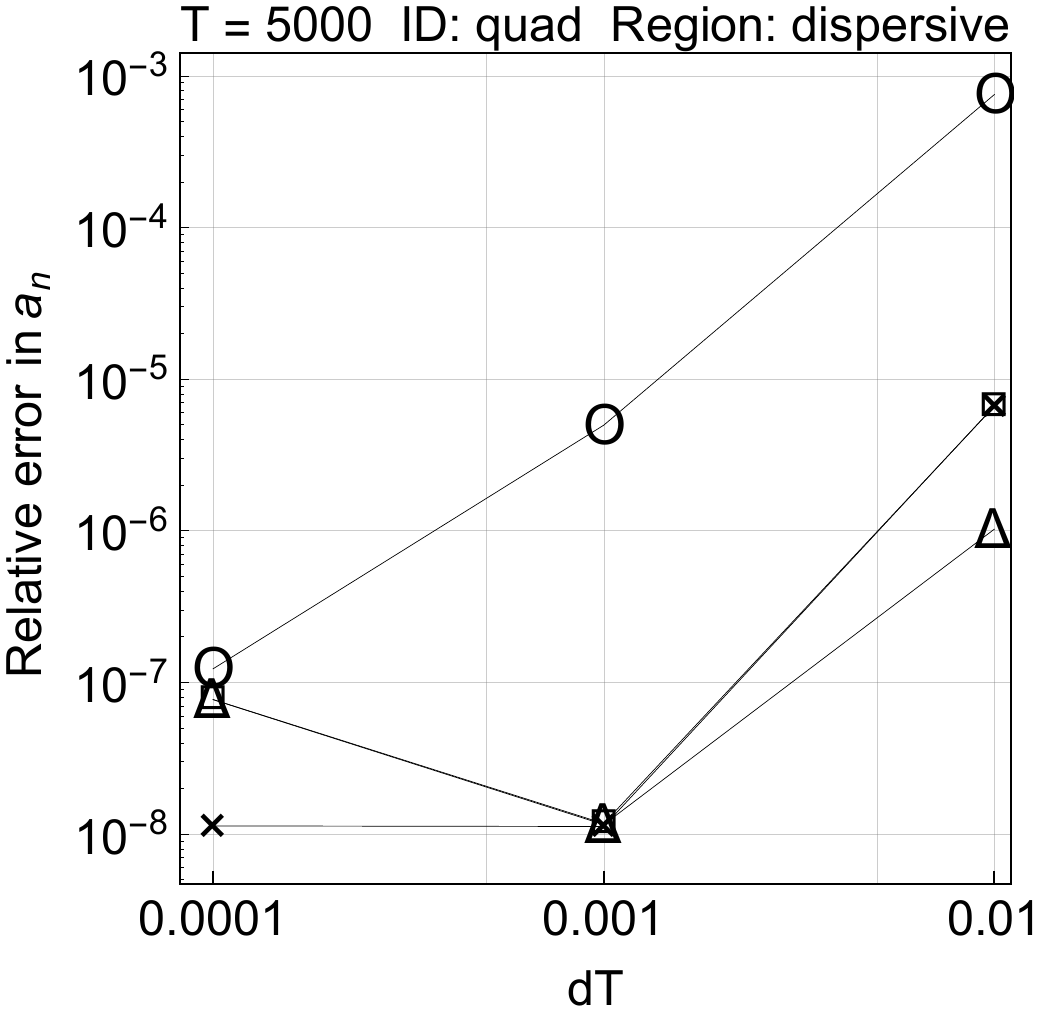}
  \includegraphics[width=.32\linewidth]{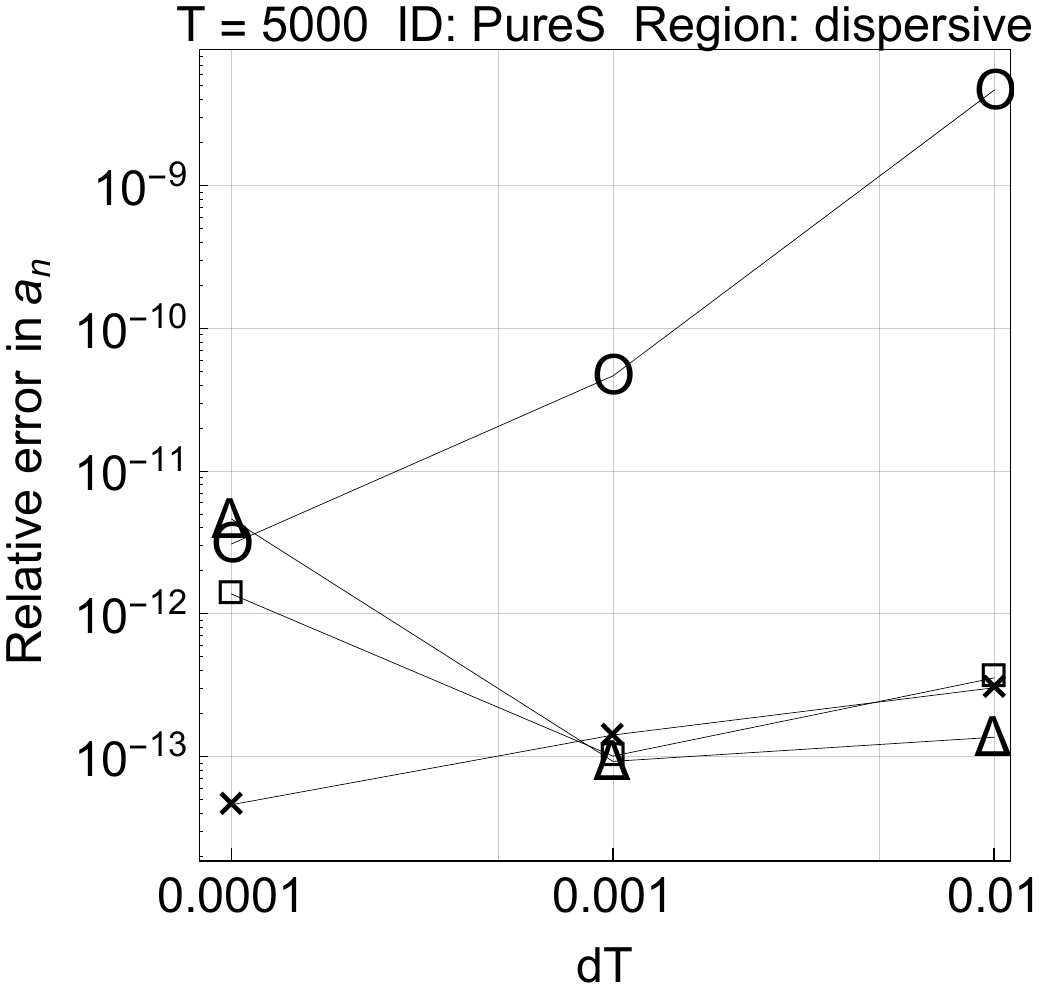}
  \caption{Relative errors in the dispersive region for the fourth-order time-stepping methods ({\tt rk4 ($\square$), rk4qp ($\times$), rkf45 ($\Delta$), ab4 ($\ocircle$)}) at $T = 5000$ plotted versus $dT$ for three choices of time step for three different choices of ID.}
    \label{f:4-disp-5000}
\end{figure}

In each of the panels of Figures~\ref{f:4-disp-2000} and \ref{f:4-disp-5000} we plot the relative error of the computed approximation of $a_n(T)$ plotted versus $dT$ for fourth-order method in the dispersive region.  We can again see that we are operating near the maximum accuracy.  We see that {\tt ab4} under-performs but not as severely as in the soliton region and {\tt rk4qp} is still almost always the method of choice in this region, at least for small time steps.  Again, the {\tt PureS} ID gives smaller errors (recall we can only measure the absolute error for {\tt PureS} ID in the dispersive region) than the other choices of ID, illustrating the importance of being able to compute these solutions accurately with the numerical IST method.

\section{Conclusions}

We have used the numerical IST method outlined in \cite{Bilman2017} to benchmark classical time-stepping routines on the Toda lattice.  It is reasonable, especially following the perturbation work in \cite{Bilman2016}, that such benchmarks will hold for other lattice equations --- not just the Toda lattice and not just integrable lattices.  The method {\tt rk4qp} appears to be the method of choices while the symplectic method {\tt sv2symp} performs surprisingly well given that it is only second-order accurate.    

We have illustrated that having a wider class of methods allows one to detect deeper differences in methods, and analyze the accumulation of round-off error while getting a handle on the maximum accuracy of a method.  We have also shown that benchmarking a method on a pure soliton initial condition can lead one to overestimate the maximum accuracy of the method.  As we consider nonlinear lattices, both the lattice equation and the choice of ID matter in the performance of a given method.

In Appendix~\ref{s:tables} we give a complete listing of relative errors for both $a_n$ and $b_n$.  In \cite{BitBenchmark} we have made available our data for the reference solutions computed with the numerical IST to allow others to benchmark their time-stepping routines.

\section*{Acknowledgments}
D.B.\ is thankful for partial support from the National Science Foundation through the grant NSF-DMS-1513054.

\appendix
\section{Butcher Tableau for Fehlberg's RKF45}\label{s:Butcher}
We provide the Butcher tableau for the constants used in the RKF45 method in this paper~\cite{Fehlberg1970}.
\begin{equation*}
\renewcommand\arraystretch{1.5}
\begin{array}{c|cccccc}
j & 1 & 2 & 3& 4& 5& 6\\
\hline
b_{2j}&\frac{1}{4} \\
b_{3j}&\frac{3}{32} & \frac{9}{32}\\
b_{4j}&\frac{1932}{2197} & \frac{-7200}{2197} & \frac{7296}{2197}\\
b_{5j}&\frac{439}{216} & -8& \frac{3680}{513} & \frac{-845}{4104} \\
b_{6j}&\frac{-8}{27} & 2 & \frac{-3544}{2565} & \frac{1859}{4104} & \frac{-11}{40}\\
\hline
c_j & \frac{25}{216} &0 & \frac{1408}{2565} &\frac{2197}{4104} &\frac{-1}{5} & 0
\end{array}
\end{equation*}

\renewcommand{\text}{\texttt}

\begin{landscape}

\section{Data for the Methods}\label{s:tables}

In the tables below, we display the data for the methods discussed in the body of the paper.  In each table we display the performance of 3-4 methods, run until $T = 1000,2000,3000$  each with three different time steps.  Each table is associated to a choice of initial data (\texttt{dirac,double,NoS,PureS,quad}) and a region (soliton region = \texttt{Sol.}, dispersive region = \texttt{Disp.}).

  \subsection{Second order --- errors for $a_n(t)$}
  \Small
  $$
  \begin{array}{|c|c|c|c|c|c|c|c|c|c|} \hline
 \text{quad/Sol.} & \text{} & \text{T = 1000} & \text{} & \text{} & \
\text{T = 2000} & \text{} & \text{} & \text{T = 5000} & \text{} \\ \hline
 \text{dT} & 0.01 & 0.001 & 0.0001 & 0.01 & 0.001 & 0.0001 & 0.01 & 0.001 & \
0.0001 \\ \hline
 \text{sv2symp} & 4.499\times 10^{-3} & 4.499\times 10^{-5} & 4.499\times \
10^{-7} & 8.914\times 10^{-3} & 8.912\times 10^{-5} & 8.912\times 10^{-7} & \
2.195\times 10^{-2} & 2.196\times 10^{-4} & 2.197\times 10^{-6} \\ \hline
 \text{midpoint} & 1.834\times 10^{-2} & 1.441\times 10^{-4} & 1.402\times \
10^{-6} & 4.222\times 10^{-2} & 2.913\times 10^{-4} & 2.782\times 10^{-6} & \
1.463\times 10^{-1} & 7.608\times 10^{-4} & 6.898\times 10^{-6} \\ \hline
 \text{midpointqp} & 2.004\times 10^{-2} & 1.561\times 10^{-4} & 1.516\times \
10^{-6} & 4.608\times 10^{-2} & 3.153\times 10^{-4} & 3.008\times 10^{-6} & \
1.593\times 10^{-1} & 8.234\times 10^{-4} & 7.456\times 10^{-6} \\ \hline
\end{array}
$$ $$
\begin{array}{|c|c|c|c|c|c|c|c|c|c|} \hline
 \text{quad/Disp.} & \text{} & \text{T = 1000} & \text{} & \text{} & \text{T \
= 2000} & \text{} & \text{} & \text{T = 5000} & \text{} \\ \hline
 \text{dT} & 0.01 & 0.001 & 0.0001 & 0.01 & 0.001 & 0.0001 & 0.01 & 0.001 & \
0.0001 \\ \hline
 \text{sv2symp} & 2.218\times 10^{-2} & 2.217\times 10^{-4} & 2.218\times \
10^{-6} & 4.43\times 10^{-2} & 4.429\times 10^{-4} & 4.428\times 10^{-6} & \
1.084\times 10^{-1} & 1.084\times 10^{-3} & 1.084\times 10^{-5} \\ \hline
 \text{midpoint} & 8.916\times 10^{-2} & 8.891\times 10^{-4} & 8.891\times \
10^{-6} & 1.781\times 10^{-1} & 1.776\times 10^{-3} & 1.775\times 10^{-5} & \
4.347\times 10^{-1} & 4.341\times 10^{-3} & 4.333\times 10^{-5} \\ \hline
 \text{midpointqp} & 8.954\times 10^{-2} & 8.896\times 10^{-4} & 8.892\times \
10^{-6} & 1.788\times 10^{-1} & 1.777\times 10^{-3} & 1.776\times 10^{-5} & \
4.364\times 10^{-1} & 4.343\times 10^{-3} & 4.34\times 10^{-5} \\ \hline
\end{array}
$$ $$
\begin{array}{|c|c|c|c|c|c|c|c|c|c|} \hline
 \text{dirac/Sol.} & \text{} & \text{T = 1000} & \text{} & \text{} & \
\text{T = 2000} & \text{} & \text{} & \text{T = 5000} & \text{} \\ \hline
 \text{dT} & 0.01 & 0.001 & 0.0001 & 0.01 & 0.001 & 0.0001 & 0.01 & 0.001 & \
0.0001 \\ \hline
 \text{sv2symp} & 2.805\times 10^{-1} & 2.937\times 10^{-3} & 2.938\times \
10^{-5} & 3.859\times 10^{-1} & 2.729\times 10^{-3} & 2.715\times 10^{-5} & \
6.994\times 10^{-1} & 7.258\times 10^{-3} & 7.347\times 10^{-5} \\ \hline
 \text{midpoint} & 1.622 & 2.808\times 10^{-2} & 9.098\times 10^{-5} & 1.451 \
& 4.956\times 10^{-2} & 1.036\times 10^{-4} & 1.428 & 1.734\times 10^{-1} & \
4.393\times 10^{-4} \\ \hline
 \text{midpointqp} & 1.634 & 3.753\times 10^{-2} & 1.285\times 10^{-4} & \
1.322 & 6.718\times 10^{-2} & 1.443\times 10^{-4} & 1.43 & 2.261\times \
10^{-1} & 5.969\times 10^{-4} \\ \hline
\end{array}
$$ $$
\begin{array}{|c|c|c|c|c|c|c|c|c|c|} \hline
 \text{dirac/Disp.} & \text{} & \text{T = 1000} & \text{} & \text{} & \text{T \
= 2000} & \text{} & \text{} & \text{T = 5000} & \text{} \\ \hline
 \text{dT} & 0.01 & 0.001 & 0.0001 & 0.01 & 0.001 & 0.0001 & 0.01 & 0.001 & \
0.0001 \\ \hline
 \text{sv2symp} & 2.229\times 10^{-2} & 2.23\times 10^{-4} & 2.23\times \
10^{-6} & 4.381\times 10^{-2} & 4.383\times 10^{-4} & 4.383\times 10^{-6} & \
1.091\times 10^{-1} & 1.089\times 10^{-3} & 1.089\times 10^{-5} \\ \hline
 \text{midpoint} & 1.044\times 10^{-1} & 9.254\times 10^{-4} & 9.138\times \
10^{-6} & 2.04\times 10^{-1} & 1.801\times 10^{-3} & 1.778\times 10^{-5} & \
5.126\times 10^{-1} & 4.444\times 10^{-3} & 4.386\times 10^{-5} \\ \hline
 \text{midpointqp} & 1.137\times 10^{-1} & 9.241\times 10^{-4} & 9.033\times \
10^{-6} & 2.242\times 10^{-1} & 1.809\times 10^{-3} & 1.768\times 10^{-5} & \
5.709\times 10^{-1} & 4.481\times 10^{-3} & 4.378\times 10^{-5} \\ \hline
\end{array}
$$ $$
\begin{array}{|c|c|c|c|c|c|c|c|c|c|} \hline
 \text{double/Sol.} & \text{} & \text{T = 1000} & \text{} & \text{} & \
\text{T = 2000} & \text{} & \text{} & \text{T = 5000} & \text{} \\ \hline
 \text{dT} & 0.01 & 0.001 & 0.0001 & 0.01 & 0.001 & 0.0001 & 0.01 & 0.001 & \
0.0001 \\ \hline
 \text{sv2symp} & 1.724\times 10^{-3} & 1.723\times 10^{-5} & 1.723\times \
10^{-7} & 3.388\times 10^{-3} & 3.387\times 10^{-5} & 3.387\times 10^{-7} & \
8.378\times 10^{-3} & 8.377\times 10^{-5} & 8.378\times 10^{-7} \\ \hline
 \text{midpoint} & 7.759\times 10^{-3} & 7.098\times 10^{-5} & 7.03\times \
10^{-7} & 1.588\times 10^{-2} & 1.404\times 10^{-4} & 1.385\times 10^{-6} & \
4.366\times 10^{-2} & 3.52\times 10^{-4} & 3.434\times 10^{-6} \\ \hline
 \text{midpointqp} & 8.85\times 10^{-3} & 8.088\times 10^{-5} & 8.013\times \
10^{-7} & 1.804\times 10^{-2} & 1.598\times 10^{-4} & 1.577\times 10^{-6} & \
4.919\times 10^{-2} & 3.999\times 10^{-4} & 3.907\times 10^{-6} \\ \hline
\end{array}
$$ $$
\begin{array}{|c|c|c|c|c|c|c|c|c|c|} \hline
 \text{double/Disp.} & \text{} & \text{T = 1000} & \text{} & \text{} & \
\text{T = 2000} & \text{} & \text{} & \text{T = 5000} & \text{} \\ \hline
 \text{dT} & 0.01 & 0.001 & 0.0001 & 0.01 & 0.001 & 0.0001 & 0.01 & 0.001 & \
0.0001 \\ \hline
 \text{sv2symp} & 2.193\times 10^{-2} & 2.193\times 10^{-4} & 2.193\times \
10^{-6} & 4.328\times 10^{-2} & 4.329\times 10^{-4} & 4.329\times 10^{-6} & \
1.099\times 10^{-1} & 1.099\times 10^{-3} & 1.099\times 10^{-5} \\ \hline
 \text{midpoint} & 8.748\times 10^{-2} & 8.763\times 10^{-4} & 8.765\times \
10^{-6} & 1.725\times 10^{-1} & 1.73\times 10^{-3} & 1.731\times 10^{-5} & \
4.369\times 10^{-1} & 4.393\times 10^{-3} & 4.392\times 10^{-5} \\ \hline
 \text{midpointqp} & 8.768\times 10^{-2} & 8.765\times 10^{-4} & 8.765\times \
10^{-6} & 1.729\times 10^{-1} & 1.731\times 10^{-3} & 1.731\times 10^{-5} & \
4.378\times 10^{-1} & 4.394\times 10^{-3} & 4.394\times 10^{-5} \\ \hline
\end{array}
$$ $$
\begin{array}{|c|c|c|c|c|c|c|c|c|c|} \hline
 \text{NoS/Sol.} & \text{} & \text{T = 1000} & \text{} & \text{} & \
\text{T = 2000} & \text{} & \text{} & \text{T = 5000} & \text{} \\ \hline
 \text{dT} & 0.01 & 0.001 & 0.0001 & 0.01 & 0.001 & 0.0001 & 0.01 & 0.001 & \
0.0001 \\ \hline
 \text{sv2symp} & 1.762\times 10^{-5} & 1.762\times 10^{-7} & 1.762\times \
10^{-9} & 2.428\times 10^{-5} & 2.427\times 10^{-7} & 2.427\times 10^{-9} & \
3.18\times 10^{-5} & 3.18\times 10^{-7} & 3.18\times 10^{-9} \\ \hline
 \text{midpoint} & 7.551\times 10^{-5} & 7.512\times 10^{-7} & 7.549\times \
10^{-9} & 1.031\times 10^{-4} & 1.026\times 10^{-6} & 1.034\times 10^{-8} & \
1.332\times 10^{-4} & 1.327\times 10^{-6} & 1.418\times 10^{-8} \\ \hline
 \text{midpointqp} & 7.82\times 10^{-5} & 7.767\times 10^{-7} & 7.762\times \
10^{-9} & 1.06\times 10^{-4} & 1.054\times 10^{-6} & 1.053\times 10^{-8} & \
1.358\times 10^{-4} & 1.352\times 10^{-6} & 1.352\times 10^{-8} \\ \hline
\end{array}
$$ $$
\begin{array}{|c|c|c|c|c|c|c|c|c|c|} \hline
 \text{NoS/Disp.} & \text{} & \text{T = 1000} & \text{} & \text{} & \text{T = \
2000} & \text{} & \text{} & \text{T = 5000} & \text{} \\ \hline
 \text{dT} & 0.01 & 0.001 & 0.0001 & 0.01 & 0.001 & 0.0001 & 0.01 & 0.001 & \
0.0001 \\ \hline
 \text{sv2symp} & 2.167\times 10^{-2} & 2.167\times 10^{-4} & 2.167\times \
10^{-6} & 4.304\times 10^{-2} & 4.303\times 10^{-4} & 4.303\times 10^{-6} & \
1.077\times 10^{-1} & 1.077\times 10^{-3} & 1.077\times 10^{-5} \\ \hline
 \text{midpoint} & 8.663\times 10^{-2} & 8.666\times 10^{-4} & 8.667\times \
10^{-6} & 1.722\times 10^{-1} & 1.721\times 10^{-3} & 1.721\times 10^{-5} & \
4.285\times 10^{-1} & 4.306\times 10^{-3} & 4.307\times 10^{-5} \\ \hline
 \text{midpointqp} & 8.673\times 10^{-2} & 8.669\times 10^{-4} & 8.669\times \
10^{-6} & 1.723\times 10^{-1} & 1.721\times 10^{-3} & 1.721\times 10^{-5} & \
4.289\times 10^{-1} & 4.307\times 10^{-3} & 4.307\times 10^{-5} \\ \hline
\end{array}
$$ $$
\begin{array}{|c|c|c|c|c|c|c|c|c|c|} \hline
 \text{PureS/Sol.} & \text{} & \text{T = 1000} & \text{} & \text{} & \
\text{T = 2000} & \text{} & \text{} & \text{T = 5000} & \text{} \\ \hline
 \text{dT} & 0.01 & 0.001 & 0.0001 & 0.01 & 0.001 & 0.0001 & 0.01 & 0.001 & \
0.0001 \\ \hline
 \text{sv2symp} & 3.707\times 10^{-4} & 3.706\times 10^{-6} & 3.706\times \
10^{-8} & 7.05\times 10^{-4} & 7.05\times 10^{-6} & 7.05\times 10^{-8} & \
1.712\times 10^{-3} & 1.712\times 10^{-5} & 1.712\times 10^{-7} \\ \hline
 \text{midpoint} & 1.485\times 10^{-3} & 1.447\times 10^{-5} & 1.442\times \
10^{-7} & 2.917\times 10^{-3} & 2.756\times 10^{-5} & 2.742\times 10^{-7} & \
7.791\times 10^{-3} & 6.756\times 10^{-5} & 6.64\times 10^{-7} \\ \hline
 \text{midpointqp} & 1.507\times 10^{-3} & 1.467\times 10^{-5} & 1.463\times \
10^{-7} & 2.964\times 10^{-3} & 2.798\times 10^{-5} & 2.781\times 10^{-7} & \
7.932\times 10^{-3} & 6.865\times 10^{-5} & 6.758\times 10^{-7} \\ \hline
\end{array}
$$ $$
\begin{array}{|c|c|c|c|c|c|c|c|c|c|} \hline
 \text{PureS/Disp.} & \text{} & \text{T = 1000} & \text{} & \text{} & \text{T \
= 2000} & \text{} & \text{} & \text{T = 5000} & \text{} \\ \hline
 \text{dT} & 0.01 & 0.001 & 0.0001 & 0.01 & 0.001 & 0.0001 & 0.01 & 0.001 & \
0.0001 \\ \hline
 \text{sv2symp} & 1.295\times 10^{-5} & 1.295\times 10^{-5} & 1.295\times \
10^{-5} & 6.467\times 10^{-6} & 6.467\times 10^{-6} & 6.467\times 10^{-6} & \
2.586\times 10^{-6} & 2.586\times 10^{-6} & 2.586\times 10^{-6} \\ \hline
 \text{midpoint} & 1.295\times 10^{-5} & 1.295\times 10^{-5} & 1.295\times \
10^{-5} & 6.469\times 10^{-6} & 6.467\times 10^{-6} & 6.467\times 10^{-6} & \
2.588\times 10^{-6} & 2.586\times 10^{-6} & 2.586\times 10^{-6} \\ \hline
 \text{midpointqp} & 1.295\times 10^{-5} & 1.295\times 10^{-5} & 1.295\times \
10^{-5} & 6.469\times 10^{-6} & 6.467\times 10^{-6} & 6.467\times 10^{-6} & \
2.588\times 10^{-6} & 2.586\times 10^{-6} & 2.586\times 10^{-6} \\ \hline
\end{array}

  $$
   \subsection{Second order --- errors for $b_n(t)$}
  $$
    \begin{array}{|c|c|c|c|c|c|c|c|c|c|} \hline
 \text{quad/Sol.} & \text{} & \text{T = 1000} & \text{} & \text{} & \
\text{T = 2000} & \text{} & \text{} & \text{T = 5000} & \text{} \\ \hline
 \text{dT} & 0.01 & 0.001 & 0.0001 & 0.01 & 0.001 & 0.0001 & 0.01 & 0.001 & \
0.0001 \\ \hline
 \text{sv2symp} & 4.471\times 10^{-3} & 4.47\times 10^{-5} & 4.47\times \
10^{-7} & 8.82\times 10^{-3} & 8.819\times 10^{-5} & 8.819\times 10^{-7} & \
2.197\times 10^{-2} & 2.196\times 10^{-4} & 2.196\times 10^{-6} \\ \hline
 \text{midpoint} & 1.823\times 10^{-2} & 1.432\times 10^{-4} & 1.393\times \
10^{-6} & 4.176\times 10^{-2} & 2.882\times 10^{-4} & 2.753\times 10^{-6} & \
1.471\times 10^{-1} & 7.608\times 10^{-4} & 6.898\times 10^{-6} \\ \hline
 \text{midpointqp} & 1.992\times 10^{-2} & 1.551\times 10^{-4} & 1.507\times \
10^{-6} & 4.558\times 10^{-2} & 3.12\times 10^{-4} & 2.976\times 10^{-6} & \
1.603\times 10^{-1} & 8.234\times 10^{-4} & 7.456\times 10^{-6} \\ \hline
\end{array}
$$ $$
\begin{array}{|c|c|c|c|c|c|c|c|c|c|} \hline
 \text{quad/Disp.} & \text{} & \text{T = 1000} & \text{} & \text{} & \text{T \
= 2000} & \text{} & \text{} & \text{T = 5000} & \text{} \\ \hline
 \text{dT} & 0.01 & 0.001 & 0.0001 & 0.01 & 0.001 & 0.0001 & 0.01 & 0.001 & \
0.0001 \\ \hline
 \text{sv2symp} & 2.18\times 10^{-2} & 2.18\times 10^{-4} & 2.181\times \
10^{-6} & 4.394\times 10^{-2} & 4.392\times 10^{-4} & 4.392\times 10^{-6} & \
1.082\times 10^{-1} & 1.082\times 10^{-3} & 1.082\times 10^{-5} \\ \hline
 \text{midpoint} & 8.761\times 10^{-2} & 8.745\times 10^{-4} & 8.745\times \
10^{-6} & 1.768\times 10^{-1} & 1.761\times 10^{-3} & 1.761\times 10^{-5} & \
4.327\times 10^{-1} & 4.337\times 10^{-3} & 4.328\times 10^{-5} \\ \hline
 \text{midpointqp} & 8.798\times 10^{-2} & 8.749\times 10^{-4} & 8.746\times \
10^{-6} & 1.775\times 10^{-1} & 1.762\times 10^{-3} & 1.761\times 10^{-5} & \
4.345\times 10^{-1} & 4.338\times 10^{-3} & 4.335\times 10^{-5} \\ \hline
\end{array}
$$ $$
\begin{array}{|c|c|c|c|c|c|c|c|c|c|} \hline
 \text{dirac/Sol.} & \text{} & \text{T = 1000} & \text{} & \text{} & \
\text{T = 2000} & \text{} & \text{} & \text{T = 5000} & \text{} \\ \hline
 \text{dT} & 0.01 & 0.001 & 0.0001 & 0.01 & 0.001 & 0.0001 & 0.01 & 0.001 & \
0.0001 \\ \hline
 \text{sv2symp} & 1.257\times 10^{-1} & 1.046\times 10^{-3} & 1.045\times \
10^{-5} & 3.659\times 10^{-1} & 4.826\times 10^{-3} & 4.835\times 10^{-5} & \
9.797\times 10^{-1} & 1.175\times 10^{-2} & 1.169\times 10^{-4} \\ \hline
 \text{midpoint} & 1.356 & 1.015\times 10^{-2} & 3.231\times 10^{-5} & 1.445 \
& 7.853\times 10^{-2} & 1.846\times 10^{-4} & 1.567 & 4.909\times 10^{-1} & \
6.994\times 10^{-4} \\ \hline
 \text{midpointqp} & 1.358 & 1.365\times 10^{-2} & 4.567\times 10^{-5} & \
1.546 & 1.026\times 10^{-1} & 2.569\times 10^{-4} & 1.638 & 6.32\times \
10^{-1} & 9.508\times 10^{-4} \\ \hline
\end{array}
$$ $$
\begin{array}{|c|c|c|c|c|c|c|c|c|c|} \hline
 \text{dirac/Disp.} & \text{} & \text{T = 1000} & \text{} & \text{} & \text{T \
= 2000} & \text{} & \text{} & \text{T = 5000} & \text{} \\ \hline
 \text{dT} & 0.01 & 0.001 & 0.0001 & 0.01 & 0.001 & 0.0001 & 0.01 & 0.001 & \
0.0001 \\ \hline
 \text{sv2symp} & 2.17\times 10^{-2} & 2.17\times 10^{-4} & 2.17\times \
10^{-6} & 4.377\times 10^{-2} & 4.377\times 10^{-4} & 4.377\times 10^{-6} & \
1.073\times 10^{-1} & 1.074\times 10^{-3} & 1.074\times 10^{-5} \\ \hline
 \text{midpoint} & 1.018\times 10^{-1} & 9.005\times 10^{-4} & 8.892\times \
10^{-6} & 2.041\times 10^{-1} & 1.799\times 10^{-3} & 1.776\times 10^{-5} & \
5.004\times 10^{-1} & 4.383\times 10^{-3} & 4.326\times 10^{-5} \\ \hline
 \text{midpointqp} & 1.109\times 10^{-1} & 8.992\times 10^{-4} & 8.79\times \
10^{-6} & 2.244\times 10^{-1} & 1.806\times 10^{-3} & 1.765\times 10^{-5} & \
5.568\times 10^{-1} & 4.42\times 10^{-3} & 4.318\times 10^{-5} \\ \hline
\end{array}
$$ $$
\begin{array}{|c|c|c|c|c|c|c|c|c|c|} \hline
 \text{double/Sol.} & \text{} & \text{T = 1000} & \text{} & \text{} & \
\text{T = 2000} & \text{} & \text{} & \text{T = 5000} & \text{} \\ \hline
 \text{dT} & 0.01 & 0.001 & 0.0001 & 0.01 & 0.001 & 0.0001 & 0.01 & 0.001 & \
0.0001 \\ \hline
 \text{sv2symp} & 1.723\times 10^{-3} & 1.723\times 10^{-5} & 1.723\times \
10^{-7} & 3.387\times 10^{-3} & 3.386\times 10^{-5} & 3.386\times 10^{-7} & \
8.373\times 10^{-3} & 8.372\times 10^{-5} & 8.372\times 10^{-7} \\ \hline
 \text{midpoint} & 7.758\times 10^{-3} & 7.096\times 10^{-5} & 7.029\times \
10^{-7} & 1.588\times 10^{-2} & 1.403\times 10^{-4} & 1.385\times 10^{-6} & \
4.363\times 10^{-2} & 3.518\times 10^{-4} & 3.432\times 10^{-6} \\ \hline
 \text{midpointqp} & 8.848\times 10^{-3} & 8.086\times 10^{-5} & 8.011\times \
10^{-7} & 1.804\times 10^{-2} & 1.597\times 10^{-4} & 1.577\times 10^{-6} & \
4.916\times 10^{-2} & 3.996\times 10^{-4} & 3.905\times 10^{-6} \\ \hline
\end{array}
$$ $$
\begin{array}{|c|c|c|c|c|c|c|c|c|c|} \hline
 \text{double/Disp.} & \text{} & \text{T = 1000} & \text{} & \text{} & \
\text{T = 2000} & \text{} & \text{} & \text{T = 5000} & \text{} \\ \hline
 \text{dT} & 0.01 & 0.001 & 0.0001 & 0.01 & 0.001 & 0.0001 & 0.01 & 0.001 & \
0.0001 \\ \hline
 \text{sv2symp} & 2.125\times 10^{-2} & 2.125\times 10^{-4} & 2.125\times \
10^{-6} & 4.355\times 10^{-2} & 4.355\times 10^{-4} & 4.355\times 10^{-6} & \
1.074\times 10^{-1} & 1.074\times 10^{-3} & 1.074\times 10^{-5} \\ \hline
 \text{midpoint} & 8.478\times 10^{-2} & 8.49\times 10^{-4} & 8.492\times \
10^{-6} & 1.737\times 10^{-1} & 1.741\times 10^{-3} & 1.741\times 10^{-5} & \
4.267\times 10^{-1} & 4.294\times 10^{-3} & 4.293\times 10^{-5} \\ \hline
 \text{midpointqp} & 8.497\times 10^{-2} & 8.492\times 10^{-4} & 8.492\times \
10^{-6} & 1.741\times 10^{-1} & 1.741\times 10^{-3} & 1.741\times 10^{-5} & \
4.277\times 10^{-1} & 4.295\times 10^{-3} & 4.295\times 10^{-5} \\ \hline
\end{array}
$$ $$
\begin{array}{|c|c|c|c|c|c|c|c|c|c|} \hline
 \text{NoS/Sol.} & \text{} & \text{T = 1000} & \text{} & \text{} & \
\text{T = 2000} & \text{} & \text{} & \text{T = 5000} & \text{} \\ \hline
 \text{dT} & 0.01 & 0.001 & 0.0001 & 0.01 & 0.001 & 0.0001 & 0.01 & 0.001 & \
0.0001 \\ \hline
 \text{sv2symp} & 1.849\times 10^{-5} & 1.849\times 10^{-7} & 1.849\times \
10^{-9} & 2.528\times 10^{-5} & 2.528\times 10^{-7} & 2.529\times 10^{-9} & \
3.275\times 10^{-5} & 3.275\times 10^{-7} & 3.276\times 10^{-9} \\ \hline
 \text{midpoint} & 7.965\times 10^{-5} & 7.925\times 10^{-7} & 7.981\times \
10^{-9} & 1.076\times 10^{-4} & 1.071\times 10^{-6} & 1.079\times 10^{-8} & \
1.373\times 10^{-4} & 1.369\times 10^{-6} & 1.461\times 10^{-8} \\ \hline
 \text{midpointqp} & 8.245\times 10^{-5} & 8.193\times 10^{-7} & 8.188\times \
10^{-9} & 1.106\times 10^{-4} & 1.1\times 10^{-6} & 1.099\times 10^{-8} & \
1.4\times 10^{-4} & 1.394\times 10^{-6} & 1.394\times 10^{-8} \\ \hline
\end{array}
$$ $$
\begin{array}{|c|c|c|c|c|c|c|c|c|c|} \hline
 \text{NoS/Disp.} & \text{} & \text{T = 1000} & \text{} & \text{} & \text{T = \
2000} & \text{} & \text{} & \text{T = 5000} & \text{} \\ \hline
 \text{dT} & 0.01 & 0.001 & 0.0001 & 0.01 & 0.001 & 0.0001 & 0.01 & 0.001 & \
0.0001 \\ \hline
 \text{sv2symp} & 2.143\times 10^{-2} & 2.143\times 10^{-4} & 2.143\times \
10^{-6} & 4.282\times 10^{-2} & 4.282\times 10^{-4} & 4.282\times 10^{-6} & \
1.089\times 10^{-1} & 1.089\times 10^{-3} & 1.089\times 10^{-5} \\ \hline
 \text{midpoint} & 8.568\times 10^{-2} & 8.569\times 10^{-4} & 8.57\times \
10^{-6} & 1.711\times 10^{-1} & 1.713\times 10^{-3} & 1.713\times 10^{-5} & \
4.328\times 10^{-1} & 4.354\times 10^{-3} & 4.354\times 10^{-5} \\ \hline
 \text{midpointqp} & 8.578\times 10^{-2} & 8.572\times 10^{-4} & 8.572\times \
10^{-6} & 1.712\times 10^{-1} & 1.713\times 10^{-3} & 1.713\times 10^{-5} & \
4.332\times 10^{-1} & 4.354\times 10^{-3} & 4.354\times 10^{-5} \\ \hline
\end{array}
$$ $$
\begin{array}{|c|c|c|c|c|c|c|c|c|c|} \hline
 \text{PureS/Sol.} & \text{} & \text{T = 1000} & \text{} & \text{} & \
\text{T = 2000} & \text{} & \text{} & \text{T = 5000} & \text{} \\ \hline
 \text{dT} & 0.01 & 0.001 & 0.0001 & 0.01 & 0.001 & 0.0001 & 0.01 & 0.001 & \
0.0001 \\ \hline
 \text{sv2symp} & 3.703\times 10^{-4} & 3.703\times 10^{-6} & 3.702\times \
10^{-8} & 7.047\times 10^{-4} & 7.047\times 10^{-6} & 7.047\times 10^{-8} & \
1.712\times 10^{-3} & 1.712\times 10^{-5} & 1.712\times 10^{-7} \\ \hline
 \text{midpoint} & 1.485\times 10^{-3} & 1.447\times 10^{-5} & 1.442\times \
10^{-7} & 2.916\times 10^{-3} & 2.755\times 10^{-5} & 2.742\times 10^{-7} & \
7.789\times 10^{-3} & 6.755\times 10^{-5} & 6.639\times 10^{-7} \\ \hline
 \text{midpointqp} & 1.506\times 10^{-3} & 1.467\times 10^{-5} & 1.463\times \
10^{-7} & 2.963\times 10^{-3} & 2.797\times 10^{-5} & 2.781\times 10^{-7} & \
7.93\times 10^{-3} & 6.864\times 10^{-5} & 6.757\times 10^{-7} \\ \hline
\end{array}
$$ $$
\begin{array}{|c|c|c|c|c|c|c|c|c|c|} \hline
 \text{PureS/Disp.} & \text{} & \text{T = 1000} & \text{} & \text{} & \text{T \
= 2000} & \text{} & \text{} & \text{T = 5000} & \text{} \\ \hline
 \text{dT} & 0.01 & 0.001 & 0.0001 & 0.01 & 0.001 & 0.0001 & 0.01 & 0.001 & \
0.0001 \\ \hline
 \text{sv2symp} & 2.589\times 10^{-5} & 2.589\times 10^{-5} & 2.589\times \
10^{-5} & 1.293\times 10^{-5} & 1.293\times 10^{-5} & 1.293\times 10^{-5} & \
5.173\times 10^{-6} & 5.172\times 10^{-6} & 5.172\times 10^{-6} \\ \hline
 \text{midpoint} & 2.59\times 10^{-5} & 2.589\times 10^{-5} & 2.589\times \
10^{-5} & 1.294\times 10^{-5} & 1.293\times 10^{-5} & 1.293\times 10^{-5} & \
5.175\times 10^{-6} & 5.172\times 10^{-6} & 5.172\times 10^{-6} \\ \hline
 \text{midpointqp} & 2.59\times 10^{-5} & 2.589\times 10^{-5} & 2.589\times \
10^{-5} & 1.294\times 10^{-5} & 1.293\times 10^{-5} & 1.293\times 10^{-5} & \
5.175\times 10^{-6} & 5.172\times 10^{-6} & 5.172\times 10^{-6} \\ \hline
\end{array}

  $$


  \subsection{Fourth order --- errors for $a_n(t)$}
  $$
  \begin{array}{|c|c|c|c|c|c|c|c|c|c|} \hline
 \text{quad/Sol.} & \text{} & \text{T = 1000} & \text{} & \text{} & \
\text{T = 2000} & \text{} & \text{} & \text{T = 5000} & \text{} \\ \hline
 \text{dT} & 0.01 & 0.001 & 0.0001 & 0.01 & 0.001 & 0.0001 & 0.01 & 0.001 & \
0.0001 \\ \hline
 \text{rk4} & 2.865\times 10^{-7} & 5.51\times 10^{-10} & 5.898\times \
10^{-10} & 6.509\times 10^{-7} & 1.476\times 10^{-10} & 2.968\times 10^{-9} \
& 2.228\times 10^{-6} & 9.889\times 10^{-10} & 2.031\times 10^{-8} \\ \hline
 \text{rk4qp} & 3.27\times 10^{-7} & 5.507\times 10^{-10} & 5.501\times \
10^{-10} & 7.39\times 10^{-7} & 5.353\times 10^{-11} & 2.099\times 10^{-11} \
& 2.511\times 10^{-6} & 1.453\times 10^{-10} & 6.118\times 10^{-11} \\ \hline
 \text{rkf45} & 6.78\times 10^{-8} & 5.502\times 10^{-10} & 5.904\times \
10^{-10} & 1.685\times 10^{-7} & 1.296\times 10^{-10} & 2.972\times 10^{-9} \
& 6.68\times 10^{-7} & 9.647\times 10^{-10} & 2.032\times 10^{-8} \\ \hline
 \text{ab4} & 3.533\times 10^{-3} & 3.53\times 10^{-5} & 3.527\times \
10^{-7} & 7.064\times 10^{-3} & 7.063\times 10^{-5} & 7.055\times 10^{-7} & \
1.745\times 10^{-2} & 1.751\times 10^{-4} & 1.749\times 10^{-6} \\ \hline
\end{array}
$$ $$
\begin{array}{|c|c|c|c|c|c|c|c|c|c|} \hline
 \text{quad/Disp.} & \text{} & \text{T = 1000} & \text{} & \text{} & \text{T \
= 2000} & \text{} & \text{} & \text{T = 5000} & \text{} \\ \hline
 \text{dT} & 0.01 & 0.001 & 0.0001 & 0.01 & 0.001 & 0.0001 & 0.01 & 0.001 & \
0.0001 \\ \hline
 \text{rk4} & 1.371\times 10^{-6} & 3.443\times 10^{-9} & 3.466\times \
10^{-9} & 2.732\times 10^{-6} & 7.455\times 10^{-9} & 8.864\times 10^{-9} & \
6.599\times 10^{-6} & 1.172\times 10^{-8} & 7.7\times 10^{-8} \\ \hline
 \text{rk4qp} & 1.381\times 10^{-6} & 3.444\times 10^{-9} & 3.468\times \
10^{-9} & 2.75\times 10^{-6} & 7.424\times 10^{-9} & 7.356\times 10^{-9} & \
6.638\times 10^{-6} & 1.12\times 10^{-8} & 1.13\times 10^{-8} \\ \hline
 \text{rkf45} & 2.129\times 10^{-7} & 3.471\times 10^{-9} & 3.466\times \
10^{-9} & 4.216\times 10^{-7} & 7.379\times 10^{-9} & 8.864\times 10^{-9} & \
1.023\times 10^{-6} & 1.152\times 10^{-8} & 7.685\times 10^{-8} \\ \hline
 \text{ab4} & 4.195\times 10^{-4} & 3.768\times 10^{-6} & 3.747\times \
10^{-8} & 5.5\times 10^{-4} & 4.546\times 10^{-6} & 4.976\times 10^{-8} & \
7.508\times 10^{-4} & 4.96\times 10^{-6} & 1.228\times 10^{-7} \\ \hline
\end{array}
$$ $$
\begin{array}{|c|c|c|c|c|c|c|c|c|c|} \hline
 \text{dirac/Sol.} & \text{} & \text{T = 1000} & \text{} & \text{} & \
\text{T = 2000} & \text{} & \text{} & \text{T = 5000} & \text{} \\ \hline
 \text{dT} & 0.01 & 0.001 & 0.0001 & 0.01 & 0.001 & 0.0001 & 0.01 & 0.001 & \
0.0001 \\ \hline
 \text{rk4} & 3.284\times 10^{-3} & 4.11\times 10^{-8} & 3.874\times \
10^{-11} & 5.916\times 10^{-3} & 6.747\times 10^{-8} & 1.03\times 10^{-10} & \
4.282\times 10^{-2} & 4.228\times 10^{-7} & 6.462\times 10^{-10} \\ \hline
 \text{rk4qp} & 4.241\times 10^{-3} & 5.391\times 10^{-8} & 2.022\times \
10^{-11} & 7.602\times 10^{-3} & 8.783\times 10^{-8} & 7.222\times 10^{-11} \
& 5.61\times 10^{-2} & 5.468\times 10^{-7} & 6.104\times 10^{-11} \\ \hline
 \text{rkf45} & 1.621\times 10^{-3} & 1.759\times 10^{-8} & 7.666\times \
10^{-11} & 3.\times 10^{-3} & 3.113\times 10^{-8} & 6.779\times 10^{-11} & \
1.947\times 10^{-2} & 2.051\times 10^{-7} & 3.397\times 10^{-10} \\ \hline
 \text{ab4} & 2.952\times 10^{-1} & 1.09\times 10^{-3} & 1.095\times \
10^{-5} & 6.553\times 10^{-1} & 1.01\times 10^{-3} & 1.016\times 10^{-5} & \
1.28 & 2.692\times 10^{-3} & 2.755\times 10^{-5} \\ \hline
\end{array}
$$ $$
\begin{array}{|c|c|c|c|c|c|c|c|c|c|} \hline
 \text{dirac/Disp.} & \text{} & \text{T = 1000} & \text{} & \text{} & \text{T \
= 2000} & \text{} & \text{} & \text{T = 5000} & \text{} \\ \hline
 \text{dT} & 0.01 & 0.001 & 0.0001 & 0.01 & 0.001 & 0.0001 & 0.01 & 0.001 & \
0.0001 \\ \hline
 \text{rk4} & 3.676\times 10^{-6} & 1.685\times 10^{-10} & 2.575\times \
10^{-10} & 7.051\times 10^{-6} & 3.221\times 10^{-10} & 8.182\times 10^{-10} \
& 1.722\times 10^{-5} & 1.122\times 10^{-9} & 6.378\times 10^{-9} \\ \hline
 \text{rk4qp} & 5.037\times 10^{-6} & 1.86\times 10^{-10} & 6.011\times \
10^{-11} & 9.665\times 10^{-6} & 3.503\times 10^{-10} & 8.565\times 10^{-11} \
& 2.362\times 10^{-5} & 8.024\times 10^{-10} & 1.398\times 10^{-10} \\ \hline
 \text{rkf45} & 1.503\times 10^{-6} & 3.749\times 10^{-11} & 2.449\times \
10^{-10} & 2.94\times 10^{-6} & 3.258\times 10^{-11} & 9.732\times 10^{-10} \
& 7.283\times 10^{-6} & 3.389\times 10^{-10} & 6.417\times 10^{-9} \\ \hline
 \text{ab4} & 5.565\times 10^{-3} & 5.217\times 10^{-5} & 5.223\times \
10^{-7} & 6.307\times 10^{-3} & 5.578\times 10^{-5} & 5.59\times 10^{-7} & \
8.02\times 10^{-3} & 6.146\times 10^{-5} & 6.204\times 10^{-7} \\ \hline
\end{array}
$$ $$
\begin{array}{|c|c|c|c|c|c|c|c|c|c|} \hline
 \text{double/Sol.} & \text{} & \text{T = 1000} & \text{} & \text{} & \
\text{T = 2000} & \text{} & \text{} & \text{T = 5000} & \text{} \\ \hline
 \text{dT} & 0.01 & 0.001 & 0.0001 & 0.01 & 0.001 & 0.0001 & 0.01 & 0.001 & \
0.0001 \\ \hline
 \text{rk4} & 1.001\times 10^{-7} & 2.442\times 10^{-11} & 1.786\times \
10^{-10} & 2.006\times 10^{-7} & 6.527\times 10^{-11} & 5.316\times 10^{-10} \
& 5.291\times 10^{-7} & 4.041\times 10^{-10} & 1.666\times 10^{-8} \\ \hline
 \text{rk4qp} & 1.129\times 10^{-7} & 9.596\times 10^{-12} & 2.287\times \
10^{-12} & 2.256\times 10^{-7} & 1.868\times 10^{-11} & 6.062\times 10^{-12} \
& 5.923\times 10^{-7} & 4.655\times 10^{-11} & 7.707\times 10^{-11} \\ \hline
 \text{rkf45} & 1.867\times 10^{-8} & 1.009\times 10^{-11} & 1.788\times \
10^{-10} & 3.843\times 10^{-8} & 3.121\times 10^{-11} & 5.324\times 10^{-10} \
& 1.083\times 10^{-7} & 3.049\times 10^{-10} & 1.666\times 10^{-8} \\ \hline
 \text{ab4} & 4.271\times 10^{-3} & 4.284\times 10^{-5} & 4.284\times \
10^{-7} & 8.436\times 10^{-3} & 8.463\times 10^{-5} & 8.461\times 10^{-7} & \
2.093\times 10^{-2} & 2.1\times 10^{-4} & 2.099\times 10^{-6} \\ \hline
\end{array}
$$ $$
\begin{array}{|c|c|c|c|c|c|c|c|c|c|} \hline
 \text{double/Disp.} & \text{} & \text{T = 1000} & \text{} & \text{} & \
\text{T = 2000} & \text{} & \text{} & \text{T = 5000} & \text{} \\ \hline
 \text{dT} & 0.01 & 0.001 & 0.0001 & 0.01 & 0.001 & 0.0001 & 0.01 & 0.001 & \
0.0001 \\ \hline
 \text{rk4} & 1.317\times 10^{-6} & 1.381\times 10^{-10} & 7.958\times \
10^{-11} & 2.595\times 10^{-6} & 2.79\times 10^{-10} & 2.274\times 10^{-10} \
& 6.58\times 10^{-6} & 1.\times 10^{-9} & 1.58\times 10^{-8} \\ \hline
 \text{rk4qp} & 1.32\times 10^{-6} & 1.312\times 10^{-10} & 1.007\times \
10^{-11} & 2.602\times 10^{-6} & 2.608\times 10^{-10} & 1.522\times 10^{-11} \
& 6.598\times 10^{-6} & 6.569\times 10^{-10} & 5.339\times 10^{-11} \\ \hline
 \text{rkf45} & 2.02\times 10^{-7} & 1.81\times 10^{-11} & 8.081\times \
10^{-11} & 3.981\times 10^{-7} & 2.337\times 10^{-11} & 2.324\times 10^{-10} \
& 1.009\times 10^{-6} & 2.204\times 10^{-10} & 1.581\times 10^{-8} \\ \hline
 \text{ab4} & 3.715\times 10^{-4} & 3.203\times 10^{-6} & 3.206\times \
10^{-8} & 4.441\times 10^{-4} & 3.416\times 10^{-6} & 3.429\times 10^{-8} & \
6.49\times 10^{-4} & 3.844\times 10^{-6} & 5.351\times 10^{-8} \\ \hline
\end{array}
$$ $$
\begin{array}{|c|c|c|c|c|c|c|c|c|c|} \hline
 \text{NoS/Sol.} & \text{} & \text{T = 1000} & \text{} & \text{} & \
\text{T = 2000} & \text{} & \text{} & \text{T = 5000} & \text{} \\ \hline
 \text{dT} & 0.01 & 0.001 & 0.0001 & 0.01 & 0.001 & 0.0001 & 0.01 & 0.001 & \
0.0001 \\ \hline
 \text{rk4} & 2.405\times 10^{-10} & 8.695\times 10^{-11} & 8.286\times \
10^{-10} & 1.726\times 10^{-10} & 1.418\times 10^{-10} & 1.287\times 10^{-9} \
& 1.282\times 10^{-10} & 4.785\times 10^{-10} & 4.865\times 10^{-9} \\ \hline
 \text{rk4qp} & 2.577\times 10^{-10} & 3.684\times 10^{-11} & 3.684\times \
10^{-11} & 1.815\times 10^{-10} & 1.949\times 10^{-11} & 1.949\times \
10^{-11} & 1.214\times 10^{-10} & 9.323\times 10^{-12} & 9.33\times 10^{-12} \
\\ \hline
 \text{rkf45} & 5.713\times 10^{-11} & 8.667\times 10^{-11} & \
8.286\times 10^{-10} & 4.106\times 10^{-11} & 1.415\times 10^{-10} & \
1.287\times 10^{-9} & 5.762\times 10^{-11} & 4.79\times 10^{-10} & \
4.865\times 10^{-9} \\ \hline
 \text{ab4} & 1.217\times 10^{-5} & 1.211\times 10^{-7} & 1.46\times \
10^{-9} & 1.149\times 10^{-5} & 1.143\times 10^{-7} & 1.727\times 10^{-9} & \
1.\times 10^{-5} & 9.939\times 10^{-8} & 4.972\times 10^{-9} \\ \hline
\end{array}
$$ $$
\begin{array}{|c|c|c|c|c|c|c|c|c|c|} \hline
 \text{NoS/Disp.} & \text{} & \text{T = 1000} & \text{} & \text{} & \text{T = \
2000} & \text{} & \text{} & \text{T = 5000} & \text{} \\ \hline
 \text{dT} & 0.01 & 0.001 & 0.0001 & 0.01 & 0.001 & 0.0001 & 0.01 & 0.001 & \
0.0001 \\ \hline
 \text{rk4} & 1.299\times 10^{-6} & 1.335\times 10^{-10} & 4.173\times \
10^{-11} & 2.578\times 10^{-6} & 2.697\times 10^{-10} & 5.083\times 10^{-11} \
& 6.454\times 10^{-6} & 6.209\times 10^{-10} & 2.603\times 10^{-10} \\ \hline
 \text{rk4qp} & 1.3\times 10^{-6} & 1.294\times 10^{-10} & 1.311\times \
10^{-11} & 2.582\times 10^{-6} & 2.607\times 10^{-10} & 1.954\times 10^{-11} \
& 6.462\times 10^{-6} & 6.337\times 10^{-10} & 4.847\times 10^{-11} \\ \hline
 \text{rkf45} & 1.997\times 10^{-7} & 2.053\times 10^{-11} & 4.3\times \
10^{-11} & 3.963\times 10^{-7} & 3.68\times 10^{-11} & 5.154\times 10^{-11} \
& 9.92\times 10^{-7} & 1.46\times 10^{-10} & 2.212\times 10^{-10} \\ \hline
 \text{ab4} & 9.795\times 10^{-5} & 6.17\times 10^{-7} & 6.166\times \
10^{-9} & 1.503\times 10^{-4} & 6.476\times 10^{-7} & 6.439\times 10^{-9} & \
3.12\times 10^{-4} & 6.84\times 10^{-7} & 6.546\times 10^{-9} \\ \hline
\end{array}
$$ $$
\begin{array}{|c|c|c|c|c|c|c|c|c|c|} \hline
 \text{PureS/Sol.} & \text{} & \text{T = 1000} & \text{} & \text{} & \
\text{T = 2000} & \text{} & \text{} & \text{T = 5000} & \text{} \\ \hline
 \text{dT} & 0.01 & 0.001 & 0.0001 & 0.01 & 0.001 & 0.0001 & 0.01 & 0.001 & \
0.0001 \\ \hline
 \text{rk4} & 4.753\times 10^{-9} & 1.185\times 10^{-10} & 2.224\times \
10^{-10} & 9.438\times 10^{-9} & 2.965\times 10^{-10} & 2.793\times 10^{-10} \
& 2.631\times 10^{-8} & 2.683\times 10^{-10} & 1.294\times 10^{-9} \\ \hline
 \text{rk4qp} & 4.881\times 10^{-9} & 2.739\times 10^{-11} & 2.744\times \
10^{-11} & 9.722\times 10^{-9} & 1.063\times 10^{-11} & 1.104\times 10^{-11} \
& 2.714\times 10^{-8} & 1.619\times 10^{-12} & 3.488\times 10^{-12} \\ \hline
 \text{rkf45} & 7.944\times 10^{-10} & 3.156\times 10^{-11} & \
8.664\times 10^{-10} & 1.693\times 10^{-9} & 3.015\times 10^{-11} & \
2.611\times 10^{-9} & 5.422\times 10^{-9} & 1.172\times 10^{-10} & \
1.265\times 10^{-9} \\ \hline
 \text{ab4} & 7.725\times 10^{-5} & 7.743\times 10^{-7} & 7.677\times \
10^{-9} & 1.579\times 10^{-4} & 1.583\times 10^{-6} & 1.592\times 10^{-8} & \
4.003\times 10^{-4} & 4.016\times 10^{-6} & 3.879\times 10^{-8} \\ \hline
\end{array}
$$ $$
\begin{array}{|c|c|c|c|c|c|c|c|c|c|} \hline
 \text{PureS/Disp.} & \text{} & \text{T = 1000} & \text{} & \text{} & \text{T \
= 2000} & \text{} & \text{} & \text{T = 5000} & \text{} \\ \hline
 \text{dT} & 0.01 & 0.001 & 0.0001 & 0.01 & 0.001 & 0.0001 & 0.01 & 0.001 & \
0.0001 \\ \hline
 \text{rk4} & 1.295\times 10^{-5} & 1.295\times 10^{-5} & 1.295\times \
10^{-5} & 6.467\times 10^{-6} & 6.467\times 10^{-6} & 6.467\times 10^{-6} & \
2.586\times 10^{-6} & 2.586\times 10^{-6} & 2.586\times 10^{-6} \\ \hline
 \text{rk4qp} & 1.295\times 10^{-5} & 1.295\times 10^{-5} & 1.295\times \
10^{-5} & 6.467\times 10^{-6} & 6.467\times 10^{-6} & 6.467\times 10^{-6} & \
2.586\times 10^{-6} & 2.586\times 10^{-6} & 2.586\times 10^{-6} \\ \hline
 \text{rkf45} & 1.295\times 10^{-5} & 1.295\times 10^{-5} & 1.295\times \
10^{-5} & 6.467\times 10^{-6} & 6.467\times 10^{-6} & 6.467\times 10^{-6} & \
2.586\times 10^{-6} & 2.586\times 10^{-6} & 2.586\times 10^{-6} \\ \hline
 \text{ab4} & 1.295\times 10^{-5} & 1.295\times 10^{-5} & 1.295\times \
10^{-5} & 6.467\times 10^{-6} & 6.467\times 10^{-6} & 6.467\times 10^{-6} & \
2.586\times 10^{-6} & 2.586\times 10^{-6} & 2.586\times 10^{-6} \\ \hline
\end{array}

  $$
   \subsection{Fourth order --- errors for $b_n(t)$}
  $$
    \begin{array}{|c|c|c|c|c|c|c|c|c|c|} \hline
 \text{quad/Sol.} & \text{} & \text{T = 1000} & \text{} & \text{} & \
\text{T = 2000} & \text{} & \text{} & \text{T = 5000} & \text{} \\ \hline
 \text{dT} & 0.01 & 0.001 & 0.0001 & 0.01 & 0.001 & 0.0001 & 0.01 & 0.001 & \
0.0001 \\ \hline
 \text{rk4} & 2.847\times 10^{-7} & 1.003\times 10^{-9} & 1.022\times \
10^{-9} & 6.441\times 10^{-7} & 1.373\times 10^{-10} & 2.712\times 10^{-9} & \
2.228\times 10^{-6} & 9.067\times 10^{-10} & 1.857\times 10^{-8} \\ \hline
 \text{rk4qp} & 3.25\times 10^{-7} & 1.003\times 10^{-9} & 1.002\times \
10^{-9} & 7.313\times 10^{-7} & 5.296\times 10^{-11} & 2.053\times 10^{-11} \
& 2.511\times 10^{-6} & 1.454\times 10^{-10} & 6.31\times 10^{-11} \\ \hline
 \text{rkf45} & 6.737\times 10^{-8} & 1.002\times 10^{-9} & 1.023\times \
10^{-9} & 1.667\times 10^{-7} & 1.186\times 10^{-10} & 2.718\times 10^{-9} & \
6.679\times 10^{-7} & 8.814\times 10^{-10} & 1.86\times 10^{-8} \\ \hline
 \text{ab4} & 3.512\times 10^{-3} & 3.508\times 10^{-5} & 3.505\times \
10^{-7} & 6.989\times 10^{-3} & 6.989\times 10^{-5} & 6.981\times 10^{-7} & \
1.747\times 10^{-2} & 1.751\times 10^{-4} & 1.749\times 10^{-6} \\ \hline
\end{array}
$$ $$
\begin{array}{|c|c|c|c|c|c|c|c|c|c|} \hline
 \text{quad/Disp.} & \text{} & \text{T = 1000} & \text{} & \text{} & \text{T \
= 2000} & \text{} & \text{} & \text{T = 5000} & \text{} \\ \hline
 \text{dT} & 0.01 & 0.001 & 0.0001 & 0.01 & 0.001 & 0.0001 & 0.01 & 0.001 & \
0.0001 \\ \hline
 \text{rk4} & 1.348\times 10^{-6} & 3.368\times 10^{-9} & 3.399\times \
10^{-9} & 2.706\times 10^{-6} & 7.313\times 10^{-9} & 8.297\times 10^{-9} & \
6.593\times 10^{-6} & 1.201\times 10^{-8} & 7.663\times 10^{-8} \\ \hline
 \text{rk4qp} & 1.358\times 10^{-6} & 3.369\times 10^{-9} & 3.406\times \
10^{-9} & 2.724\times 10^{-6} & 7.302\times 10^{-9} & 7.283\times 10^{-9} & \
6.632\times 10^{-6} & 1.119\times 10^{-8} & 1.117\times 10^{-8} \\ \hline
 \text{rkf45} & 2.101\times 10^{-7} & 3.413\times 10^{-9} & 3.399\times \
10^{-9} & 4.185\times 10^{-7} & 7.276\times 10^{-9} & 8.305\times 10^{-9} & \
1.022\times 10^{-6} & 1.175\times 10^{-8} & 7.646\times 10^{-8} \\ \hline
 \text{ab4} & 4.187\times 10^{-4} & 3.77\times 10^{-6} & 3.745\times \
10^{-8} & 5.455\times 10^{-4} & 4.51\times 10^{-6} & 4.83\times 10^{-8} & \
7.542\times 10^{-4} & 4.992\times 10^{-6} & 1.231\times 10^{-7} \\ \hline
\end{array}
$$ $$
\begin{array}{|c|c|c|c|c|c|c|c|c|c|} \hline
 \text{dirac/Sol.} & \text{} & \text{T = 1000} & \text{} & \text{} & \
\text{T = 2000} & \text{} & \text{} & \text{T = 5000} & \text{} \\ \hline
 \text{dT} & 0.01 & 0.001 & 0.0001 & 0.01 & 0.001 & 0.0001 & 0.01 & 0.001 & \
0.0001 \\ \hline
 \text{rk4} & 1.164\times 10^{-3} & 1.46\times 10^{-8} & 1.849\times \
10^{-11} & 1.071\times 10^{-2} & 1.202\times 10^{-7} & 1.683\times 10^{-10} \
& 6.202\times 10^{-2} & 6.725\times 10^{-7} & 9.511\times 10^{-10} \\ \hline
 \text{rk4qp} & 1.503\times 10^{-3} & 1.915\times 10^{-8} & 1.413\times \
10^{-11} & 1.382\times 10^{-2} & 1.565\times 10^{-7} & 1.287\times 10^{-10} \
& 7.915\times 10^{-2} & 8.697\times 10^{-7} & 9.272\times 10^{-11} \\ \hline
 \text{rkf45} & 5.764\times 10^{-4} & 6.248\times 10^{-9} & 3.149\times \
10^{-11} & 5.305\times 10^{-3} & 5.547\times 10^{-8} & 8.918\times 10^{-11} \
& 3.248\times 10^{-2} & 3.262\times 10^{-7} & 3.574\times 10^{-10} \\ \hline
 \text{ab4} & 9.661\times 10^{-2} & 3.876\times 10^{-4} & 3.889\times \
10^{-6} & 9.443\times 10^{-1} & 1.794\times 10^{-3} & 1.809\times 10^{-5} & \
1.495 & 4.31\times 10^{-3} & 4.382\times 10^{-5} \\ \hline
\end{array}
$$ $$
\begin{array}{|c|c|c|c|c|c|c|c|c|c|} \hline
 \text{dirac/Disp.} & \text{} & \text{T = 1000} & \text{} & \text{} & \text{T \
= 2000} & \text{} & \text{} & \text{T = 5000} & \text{} \\ \hline
 \text{dT} & 0.01 & 0.001 & 0.0001 & 0.01 & 0.001 & 0.0001 & 0.01 & 0.001 & \
0.0001 \\ \hline
 \text{rk4} & 3.579\times 10^{-6} & 1.693\times 10^{-10} & 2.559\times \
10^{-10} & 7.044\times 10^{-6} & 3.165\times 10^{-10} & 8.138\times 10^{-10} \
& 1.699\times 10^{-5} & 1.126\times 10^{-9} & 6.312\times 10^{-9} \\ \hline
 \text{rk4qp} & 4.905\times 10^{-6} & 1.876\times 10^{-10} & 6.095\times \
10^{-11} & 9.656\times 10^{-6} & 3.43\times 10^{-10} & 8.417\times 10^{-11} \
& 2.33\times 10^{-5} & 8.094\times 10^{-10} & 1.479\times 10^{-10} \\ \hline
 \text{rkf45} & 1.464\times 10^{-6} & 3.167\times 10^{-11} & 2.43\times \
10^{-10} & 2.938\times 10^{-6} & 3.728\times 10^{-11} & 9.681\times 10^{-10} \
& 7.185\times 10^{-6} & 3.512\times 10^{-10} & 6.348\times 10^{-9} \\ \hline
 \text{ab4} & 5.444\times 10^{-3} & 5.106\times 10^{-5} & 5.111\times \
10^{-7} & 6.308\times 10^{-3} & 5.581\times 10^{-5} & 5.592\times 10^{-7} & \
7.896\times 10^{-3} & 6.046\times 10^{-5} & 6.104\times 10^{-7} \\ \hline
\end{array}
$$ $$
\begin{array}{|c|c|c|c|c|c|c|c|c|c|} \hline
 \text{double/Sol.} & \text{} & \text{T = 1000} & \text{} & \text{} & \
\text{T = 2000} & \text{} & \text{} & \text{T = 5000} & \text{} \\ \hline
 \text{dT} & 0.01 & 0.001 & 0.0001 & 0.01 & 0.001 & 0.0001 & 0.01 & 0.001 & \
0.0001 \\ \hline
 \text{rk4} & 1.001\times 10^{-7} & 2.425\times 10^{-11} & 1.776\times \
10^{-10} & 2.006\times 10^{-7} & 6.516\times 10^{-11} & 5.24\times 10^{-10} \
& 5.287\times 10^{-7} & 3.88\times 10^{-10} & 1.59\times 10^{-8} \\ \hline
 \text{rk4qp} & 1.128\times 10^{-7} & 9.553\times 10^{-12} & 2.335\times \
10^{-12} & 2.255\times 10^{-7} & 1.866\times 10^{-11} & 6.009\times 10^{-12} \
& 5.919\times 10^{-7} & 4.644\times 10^{-11} & 7.897\times 10^{-11} \\ \hline
 \text{rkf45} & 1.867\times 10^{-8} & 9.875\times 10^{-12} & 1.779\times \
10^{-10} & 3.842\times 10^{-8} & 3.104\times 10^{-11} & 5.248\times 10^{-10} \
& 1.083\times 10^{-7} & 2.96\times 10^{-10} & 1.589\times 10^{-8} \\ \hline
 \text{ab4} & 4.27\times 10^{-3} & 4.283\times 10^{-5} & 4.283\times \
10^{-7} & 8.434\times 10^{-3} & 8.461\times 10^{-5} & 8.459\times 10^{-7} & \
2.091\times 10^{-2} & 2.098\times 10^{-4} & 2.098\times 10^{-6} \\ \hline
\end{array}
$$ $$
\begin{array}{|c|c|c|c|c|c|c|c|c|c|} \hline
 \text{double/Disp.} & \text{} & \text{T = 1000} & \text{} & \text{} & \
\text{T = 2000} & \text{} & \text{} & \text{T = 5000} & \text{} \\ \hline
 \text{dT} & 0.01 & 0.001 & 0.0001 & 0.01 & 0.001 & 0.0001 & 0.01 & 0.001 & \
0.0001 \\ \hline
 \text{rk4} & 1.275\times 10^{-6} & 1.357\times 10^{-10} & 7.93\times \
10^{-11} & 2.61\times 10^{-6} & 2.795\times 10^{-10} & 2.291\times 10^{-10} \
& 6.431\times 10^{-6} & 9.859\times 10^{-10} & 1.544\times 10^{-8} \\ \hline
 \text{rk4qp} & 1.278\times 10^{-6} & 1.287\times 10^{-10} & 8.07\times \
10^{-12} & 2.617\times 10^{-6} & 2.611\times 10^{-10} & 1.358\times 10^{-11} \
& 6.448\times 10^{-6} & 6.493\times 10^{-10} & 4.824\times 10^{-11} \\ \hline
 \text{rkf45} & 1.956\times 10^{-7} & 1.569\times 10^{-11} & 8.033\times \
10^{-11} & 4.003\times 10^{-7} & 2.324\times 10^{-11} & 2.334\times 10^{-10} \
& 9.861\times 10^{-7} & 2.212\times 10^{-10} & 1.546\times 10^{-8} \\ \hline
 \text{ab4} & 3.616\times 10^{-4} & 3.123\times 10^{-6} & 3.126\times \
10^{-8} & 4.477\times 10^{-4} & 3.446\times 10^{-6} & 3.459\times 10^{-8} & \
6.343\times 10^{-4} & 3.762\times 10^{-6} & 5.233\times 10^{-8} \\ \hline
\end{array}
$$ $$
\begin{array}{|c|c|c|c|c|c|c|c|c|c|} \hline
 \text{NoS/Sol.} & \text{} & \text{T = 1000} & \text{} & \text{} & \
\text{T = 2000} & \text{} & \text{} & \text{T = 5000} & \text{} \\ \hline
 \text{dT} & 0.01 & 0.001 & 0.0001 & 0.01 & 0.001 & 0.0001 & 0.01 & 0.001 & \
0.0001 \\ \hline
 \text{rk4} & 2.315\times 10^{-10} & 8.416\times 10^{-11} & 8.496\times \
10^{-10} & 1.709\times 10^{-10} & 1.501\times 10^{-10} & 1.363\times 10^{-9} \
& 1.287\times 10^{-10} & 4.976\times 10^{-10} & 5.047\times 10^{-9} \\ \hline
 \text{rk4qp} & 2.47\times 10^{-10} & 2.2\times 10^{-11} & 2.201\times \
10^{-11} & 1.786\times 10^{-10} & 1.828\times 10^{-11} & 1.828\times \
10^{-11} & 1.21\times 10^{-10} & 1.335\times 10^{-11} & 1.335\times 10^{-11} \
\\ \hline
 \text{rkf45} & 4.745\times 10^{-11} & 8.395\times 10^{-11} & \
8.496\times 10^{-10} & 4.034\times 10^{-11} & 1.497\times 10^{-10} & \
1.363\times 10^{-9} & 6.015\times 10^{-11} & 4.979\times 10^{-10} & \
5.047\times 10^{-9} \\ \hline
 \text{ab4} & 1.246\times 10^{-5} & 1.24\times 10^{-7} & 1.511\times \
10^{-9} & 1.175\times 10^{-5} & 1.168\times 10^{-7} & 1.788\times 10^{-9} & \
1.017\times 10^{-5} & 1.011\times 10^{-7} & 5.151\times 10^{-9} \\ \hline
\end{array}
$$ $$
\begin{array}{|c|c|c|c|c|c|c|c|c|c|} \hline
 \text{NoS/Disp.} & \text{} & \text{T = 1000} & \text{} & \text{} & \text{T = \
2000} & \text{} & \text{} & \text{T = 5000} & \text{} \\ \hline
 \text{dT} & 0.01 & 0.001 & 0.0001 & 0.01 & 0.001 & 0.0001 & 0.01 & 0.001 & \
0.0001 \\ \hline
 \text{rk4} & 1.282\times 10^{-6} & 1.325\times 10^{-10} & 4.284\times \
10^{-11} & 2.566\times 10^{-6} & 2.653\times 10^{-10} & 4.116\times 10^{-11} \
& 6.524\times 10^{-6} & 6.385\times 10^{-10} & 2.428\times 10^{-10} \\ \hline
 \text{rk4qp} & 1.284\times 10^{-6} & 1.278\times 10^{-10} & 9.143\times \
10^{-12} & 2.57\times 10^{-6} & 2.546\times 10^{-10} & 1.756\times 10^{-11} \
& 6.532\times 10^{-6} & 6.535\times 10^{-10} & 3.917\times 10^{-11} \\ \hline
 \text{rkf45} & 1.971\times 10^{-7} & 1.87\times 10^{-11} & 4.349\times \
10^{-11} & 3.946\times 10^{-7} & 3.901\times 10^{-11} & 4.469\times 10^{-11} \
& 1.003\times 10^{-6} & 1.359\times 10^{-10} & 2.087\times 10^{-10} \\ \hline
 \text{ab4} & 9.781\times 10^{-5} & 6.202\times 10^{-7} & 6.199\times \
10^{-9} & 1.487\times 10^{-4} & 6.416\times 10^{-7} & 6.377\times 10^{-9} & \
3.154\times 10^{-4} & 6.868\times 10^{-7} & 6.585\times 10^{-9} \\ \hline
\end{array}
$$ $$
\begin{array}{|c|c|c|c|c|c|c|c|c|c|} \hline
 \text{PureS/Sol.} & \text{} & \text{T = 1000} & \text{} & \text{} & \
\text{T = 2000} & \text{} & \text{} & \text{T = 5000} & \text{} \\ \hline
 \text{dT} & 0.01 & 0.001 & 0.0001 & 0.01 & 0.001 & 0.0001 & 0.01 & 0.001 & \
0.0001 \\ \hline
 \text{rk4} & 4.752\times 10^{-9} & 1.161\times 10^{-10} & 2.201\times \
10^{-10} & 9.436\times 10^{-9} & 2.965\times 10^{-10} & 2.74\times 10^{-10} \
& 2.63\times 10^{-8} & 2.681\times 10^{-10} & 1.292\times 10^{-9} \\ \hline
 \text{rk4qp} & 4.879\times 10^{-9} & 1.034\times 10^{-11} & 1.054\times \
10^{-11} & 9.72\times 10^{-9} & 5.901\times 10^{-12} & 7.241\times 10^{-12} \
& 2.713\times 10^{-8} & 2.156\times 10^{-12} & 1.005\times 10^{-11} \\ \hline
 \text{rkf45} & 7.937\times 10^{-10} & 1.935\times 10^{-11} & \
8.682\times 10^{-10} & 1.692\times 10^{-9} & 2.888\times 10^{-11} & \
2.609\times 10^{-9} & 5.421\times 10^{-9} & 1.172\times 10^{-10} & \
1.264\times 10^{-9} \\ \hline
 \text{ab4} & 7.724\times 10^{-5} & 7.741\times 10^{-7} & 7.675\times \
10^{-9} & 1.579\times 10^{-4} & 1.583\times 10^{-6} & 1.592\times 10^{-8} & \
4.002\times 10^{-4} & 4.015\times 10^{-6} & 3.878\times 10^{-8} \\ \hline
\end{array}
$$ $$
\begin{array}{|c|c|c|c|c|c|c|c|c|c|} \hline
 \text{PureS/Disp.} & \text{} & \text{T = 1000} & \text{} & \text{} & \text{T \
= 2000} & \text{} & \text{} & \text{T = 5000} & \text{} \\ \hline
 \text{dT} & 0.01 & 0.001 & 0.0001 & 0.01 & 0.001 & 0.0001 & 0.01 & 0.001 & \
0.0001 \\ \hline
 \text{rk4} & 2.589\times 10^{-5} & 2.589\times 10^{-5} & 2.589\times \
10^{-5} & 1.293\times 10^{-5} & 1.293\times 10^{-5} & 1.293\times 10^{-5} & \
5.172\times 10^{-6} & 5.172\times 10^{-6} & 5.172\times 10^{-6} \\ \hline
 \text{rk4qp} & 2.589\times 10^{-5} & 2.589\times 10^{-5} & 2.589\times \
10^{-5} & 1.293\times 10^{-5} & 1.293\times 10^{-5} & 1.293\times 10^{-5} & \
5.172\times 10^{-6} & 5.172\times 10^{-6} & 5.172\times 10^{-6} \\ \hline
 \text{rkf45} & 2.589\times 10^{-5} & 2.589\times 10^{-5} & 2.589\times \
10^{-5} & 1.293\times 10^{-5} & 1.293\times 10^{-5} & 1.293\times 10^{-5} & \
5.172\times 10^{-6} & 5.172\times 10^{-6} & 5.172\times 10^{-6} \\ \hline
 \text{ab4} & 2.589\times 10^{-5} & 2.589\times 10^{-5} & 2.589\times \
10^{-5} & 1.294\times 10^{-5} & 1.293\times 10^{-5} & 1.293\times 10^{-5} & \
5.173\times 10^{-6} & 5.172\times 10^{-6} & 5.172\times 10^{-6} \\ \hline
\end{array}

  $$
\end{landscape}

\bibliographystyle{siam}
\bibliography{DBTT-Benchmark}
\end{document}